\documentclass[11pt]{article}
  \usepackage{amsthm}
  \usepackage{amssymb}
  \usepackage{amsmath}
  \usepackage{eucal}
  \usepackage{mathrsfs}
  \usepackage[all]{xy}
 \usepackage{bbm}
 \usepackage{geometry}
\newcommand{\Po}{K} % poset 
\newcommand{\dc}{o} % aperto
\newcommand{\si}{\sigma}
\newcommand{\Si}{\Sigma}

\newcommand{\io}{\iota}

\newcommand{\Hom}{\mathrm{Hom}}

\newcommand{\Ccal}{\mathcal{C}}
\newcommand{\Ecal}{\mathcal{E}} % vector bundle
\newcommand{\Fcal}{\mathcal{F}} 
\newcommand{\Hcal}{\mathcal{H}} 
 
\newcommand{\Lcal}{\mathcal{L}}
\newcommand{\Pcal}{\mathcal{P}} % Fibrato principale 

\newcommand{\Tcal}{\mathcal{T}} % fibrato banale
\newcommand{\Ucal}{\mathcal{U}} % covering 
\newcommand{\Vcal}{\mathcal{V}} % 
\newcommand{\Xcal}{\mathcal{X}}

\newcommand{\dr}{\mathrm{d}} % standard fiber 
\newcommand{\Zrm}{\mathrm{Z}}
\newcommand{\Crm}{\mathrm{C}} %
\newcommand{\Hrm}{\mathrm{H}} %
\newcommand{\Brm}{\mathrm{B}} %
\newcommand{\Rrm}{\mathrm{R}} %

\newcommand{\cf}{\mathfrak{c}}

\newcommand{\Zo}{\mathbb{Z}} % Insieme numeri interi
\newcommand{\Co}{\mathbb{C}} % Insieme numeri complessi
\newcommand{\Uo}{\mathbb{U}} % gruppo unitario
\newcommand{\To}{\mathbb{T}} % toro
\newcommand{\No}{\mathbb{N}} % numeri naturali
\newcommand{\mB}{\mathfrak{B}} % Operatori limitati
\newcommand{\Ro}{\mathbb R}
\newcommand{\mcPE}{ \mathcal{PE}_\prec }   % fibrato proiettivo
\newcommand{\mPE}{ PE_\prec }               % supporto poset fibrato proiettivo
\newcommand{\mPEa}{ PE_{\prec,a} } 
\newcommand{\Krm}{{\bf{K}}}

\author{John E. Roberts$^{(1)}$, Giuseppe Ruzzi$^{(1)}$, 
        Ezio Vasselli$^{(2)}$ \\[5pt] 
\small{$^{(1)}$ Dipartimento di Matematica, Universit\`a di 
       Roma ``Tor Vergata'' }\\
      \small{Via della Ricerca Scientifica I-00133, Roma,  Italy} \\[5pt] 
\small{$^{(2)}$ Dipartimento di Matematica, Universit\`a di 
       Roma ``La Sapienza'' }\\
      \small{Piazzale Aldo Moro 5, I-00185 Roma, Italy}\\[5pt] 
    \small{\texttt{roberts@mat.uniroma2.it}} \ , \
            \small{\texttt{ruzzi@mat.uniroma2.it}} \ , \  
            \small{\texttt{ezio.vasselli@gmail.com}}}

\title{Net bundles over posets and K-theory}
\date{}

\numberwithin{equation}{section}

\begin{document}

\maketitle

\begin{abstract}
We continue studying net bundles over partially ordered sets 
(posets), defined as the analogues of ordinary fibre bundles. To this end, 
we analyze the connection between homotopy, net homology and net cohomology 
of a poset, giving versions of classical Hurewicz theorems.
Focusing our attention on Hilbert net bundles, we define the K-theory
of a poset and introduce functions over the homotopy groupoid 
satisfying the same formal properties as Chern classes.
As when the given poset is a base for the topology of a
space, our results apply to the category of locally constant bundles.

\

\noindent {\bf MSC-class}: 13D15; 05E25; 06A11.

\end{abstract}

\theoremstyle{plain}
\newtheorem{df}{Definition}[section]
\newtheorem{teo}[df]{Theorem}
\newtheorem{prop}[df]{Proposition}
\newtheorem{cor}[df]{Corollary}
\newtheorem{lemma}[df]{Lemma}

\theoremstyle{definition}
\newtheorem{oss}[df]{Remark}
\newtheorem{ex}[df]{Example}

%*************************************************************
\tableofcontents
% \markboth{Contents}{Contents}
%*************************************************************

\section{Introduction} 

This paper continues the discussion of invariants of a partially ordered set 
begun in \cite{RR06,RRV07}, comparing them with the analogous topological 
invariants when the partially ordered set is a suitable base for 
the topology of a topological space ordered under inclusion. Whereas previously, 
the emphasis has been on the fundamental groupoid, connections, curvature, 
simplicial cohomology and \v Cech cohomology, it is now on homology, locally 
constant cohomology and K-theory with the fundamental groupoid continuing its 
ubiquitous role.

Thus we prove an analogue of the Hurewicz isomorphism showing that the first 
homology group of a poset with values in $\Zo$ is isomorphic to its Abelianized 
homotopy group. We further consider an arcwise locally contractible space and a 
poset made up of a base of contractible open sets ordered under inclusion 
and show that the first homology groups with values in $\Zo$ of the topological 
space coincides with that of the poset as do the first cohomology groups with 
values in an Abelian group.\smallskip 

The pivotal notions are those of of net bundle and quasinet bundle. These are 
bundles over a poset where the fibres are objects in a category together with a 
functor from the poset mapping an element of the poset to the corresponding fibre. 
In a quasinet bundle, the functor takes values in the monomorphisms, in a net bundle 
in the isomorphisms.\smallskip

We show that a poset net bundle with pathwise connected fibres gives rise to a short 
exact sequence of homotopy groups as does a net bundle of topological spaces with 
connected and locally contractible fibres.

We define a functor  from the category of net bundles of topological spaces to the 
category of fibre bundles over the poset with the Alexandroff topology. The category 
of $1$--cocycles with values in a group of the global space of the original net bundle 
regarded as a net bundle over the poset with the opposite ordering is equivalent to 
the category of homomorphisms of the homotopy group of the fibre bundle with values 
in the group. The fibre bundle associated with a principal $G$--net bundle has locally 
constant transition functions.\smallskip

The remainder of the paper is largely devoted to studying Hilbert net bundles and 
their K-theory. Here there is a functor from Hilbert net bundles to vector 
bundles over the poset with the Alexandroff topology. The category of Hilbert net 
bundles is equivalent to the category of unitary finite-dimensional representations 
of the homotopy group of the poset. Consequently, its K-ring is isomorphic 
to the representation ring of its homotopy group. But it is also equivalent to 
the category of locally constant vector bundles over the poset with the Alexandroff 
topology. This leads to an isomorphism of the corresponding K-rings.
\smallskip 

We prove an analogue of the Thom isomorphism asserting that the first cohomology of 
the poset with values in a group is isomorphic to the first cohomology of a projective 
net bundle with values in the same group.\smallskip

The first Chern class of a Hilbert net bundle is defined as an 
element of the net cohomology group $\Hrm^1(\Po,\To)$. 
After this, the Chern K-classes are defined as elements 
of the reduced K-theory and Chern functions as complex-valued functions on 
the homotopy 
classes of paths, allowing one to recover the first Chern class by  
evaluating a suitable polynomial. When  $\Po$ is a 
base for the topology of
a space $M$, our results describe locally constant vector bundles over $M$,
the first Chern class is an element of the singular cohomology $\Hrm^1(M,\To)$
and the Chern functions are defined on the homotopy groupoid of $M$.
\smallskip

Our results imply 
that when the poset $\Po$ is a base for the topology of a 
manifold, each finite-dimensional Banach net bundle over $\Po$ has well-defined 
secondary characteristic classes (see \cite{CS85}, \cite[\S 4.19]{Kar2},\cite[\S 1(g)]{BL95}). 
It would be interesting to find a description of such classes in terms of 
the properties of $\Po$ such as its net cohomology.

Finally, we would like to mention that our approach to the homotopy group
may be regarded as a particular case of the fundamental work of Quillen on higher algebraic K-theory (\cite{Qui}),
focused on generic categories instead of posets. 
Neverthless our approach has the advantage of being expressed in a language derived
from algebraic quantum field theory, which motivated our work.

\bigskip 

An appendix is devoted to describing the simplicial sets associated with a poset, 
explaining the related notion of homotopy and proving that the homotopy group of a 
product of symmetric simplicial sets is equal to the product of the homotopy groups.
\smallskip

\section{Background and Notations.}

For the basic notions on the simplicial set of a poset $\Po$ and its homotopy
theory, including the definition of path, we refer the reader to 
\cite{Ruz05,RR06,RRV07}. However the basic definitions and terminology 
can also be found in the Appendix.

The fundamental covering of $\Po$ is given by the family of sets
\[
V_a := \{ \dc \in \Po : a \leq \dc \} 
\ \ , \ \
a \in \Sigma_0(\Po)
\ ,
\]
which provide a base for the Alexandroff topology of $\Po$. We denote
the associated space by $\tau \Po$ (see \cite[\S 2.3]{RRV07}).

In order to simplify the exposition, in the present paper {\em we shall
always assume that the poset $\Po$ is pathwise connected}. This implies
that the isomorphism class of the homotopy group $\pi_1(\Po,a)$ does not
depend on the choice of $a \in \Si_0(\Po)$.

In the same way, each space $M$ is assumed to be arcwise connected,
and its homotopy group is denoted by $\pi_1(M)$. A pivotal result that
will be used in the present paper is the following: when $M$ is Hausdorff,
and $\Po$ is a base of arcwise and simply connected open subsets of $M$
ordered under inclusion, there is an isomorphism 
$\pi_1(M) \simeq \pi_1(\Po,a)$, $a \in \Si_0(\Po)$ 
(see \cite[Thm.2.18]{Ruz05}).

We denote the identity map of a set $S$ by $id_S$.

A multiplicative semigroup of polynomials with coefficients in a ring $R$
is defined by 
\[
1 + h R [[h]]
\ := \
\{
1 + \sum_{k=1}^r x_k h^k
\ \ , \ \
r \in \No
\ , \
x_k \in R \
\}
\ .
\]

%*************************************************************

\section{Abelian (co)homology.}
\label{Ba}

The relation between homotopy and Abelian (co)homology of a poset can be deduced by the Quillen's paper \cite{Qui}. 
However, we prefer to analyze this relation in terms of a simplicial set associated with the poset, 
an approach closer to the language used in the present paper than that used by Quillen.  
The main aim of this section is to give some versions of the classical Hurewicz theorems. \\
\indent We refer the reader to the Appendix \ref{sec_ssets} for the definition of the simplicial set associated to the poset
and for the corresponding homotopy.

\subsection{Net cohomology.}
The net cohomology of the poset $\Po$ with values in an Abelian group $A$, 
written additively, 
is the cohomology of the simplicial set 
$\widetilde \Si_*(\Po)$  with values in  $A$ (see Appendix \ref{sec_ssets}). To be precise,  
one can define  the set $\Crm^n(\Po,A)$ of $n$--cochains of $\Po$ 
with values in $A$ 
as the set of functions $v: \widetilde \Si_n(K) \to A$.
$\Crm^n(\Po,A)$ inherits from $A$ the structure of an Abelian group:
\[
(v + w) (x) := v(x) + w(x)
\  , \ 
(- v)(x) := -v(x)
\ , \ 
x\in \widetilde \Si_n(\Po) \ ,  
\]
for any $v,w \in \Crm^n(\Po,A)$. The coboundary operator $\dr$ defined by 
\[
\dr v (x) =  \sum^n_{k=0} (-1)^k  \ v(\partial_k x), \qquad 
x\in \widetilde \Si_n(K),
\]
is a  mapping $\dr:\Crm^n(\Po,A)\rightarrow\Crm^{n+1}(\Po,A)$ satisfying
the equation 
$\dr\dr v = \io$, $v\in \Crm^n(\Po,A)$,
where $\io$ is the trivial cochain.  Clearly $\dr$ is a group morphism. 
An $n$--cochain $z$  
is said to be an $n$--{\em cocycle} if belongs to the kernel of $\dr$.
We denote  the group of $n$--cocycles by $\Zrm^n(\Po,A)$. Since 
$\dr \Crm^{n-1}( \Po , A )$ is a subgroup $\Brm^n(\Po,A)$ of $\Crm^n(\Po,A)$, 
we define the {\em net cohomology with coefficients in $A$} 
\[
\Hrm^n(\Po , A) \ := \ \Zrm^n(\Po,A) / \Brm^n(\Po,A) \ .
\]
The functors $\Hrm^n( - , - )$, $n \in \No$, are contravariant w.r.t 
the poset and covariant w.r.t. to the group; moreover, net cohomology 
has long exact sequences (see \cite[Lemma 1.1]{LR}).

\begin{oss}
\label{oss_h1A}
Let $\Po$ be connected, $a \in \Si_0(\Po)$ and
${\mathrm{Hom}}( \pi_1(\Po,a) , A )$ denote the set of morphisms 
from $\pi_1(\Po,a)$ to $A$. Using \cite[Prop.3.8]{RR06} and the fact 
that $A$ is Abelian, there is an isomorphism 
$\Hrm^1 (\Po,A)  \simeq  {\mathrm{Hom}}( \pi_1(\Po,a) , A )$.
This allows one to compute several net cohomology groups.
\end{oss}

\subsection{Net homology.}
\label{sec_homology}

In the present section, we introduce the notion of the net homology of a poset.

Let $A$ be an Abelian group with identity $0 \in A$.
For every $n \in \No$, we let $\Crm_n(\Po,A)$ denote the free Abelian group 
generated by formal linear combinations of elements of $\widetilde \Si_n (\Po)$
with coefficients in $A$, and define the boundary morphism as the $A$-linear map
\begin{equation}
\label{eq_bound}
{\mathrm{b}}_n : \Crm_n(\Po,A) \to \Crm_{n-1}(\Po,A)
\ \ , \ \
{\mathrm{b}}_n x 
\ := \ 
\sum_{k=0}^n (-1)^k \partial_k x 
\ .
\end{equation}
It is clear that ${\mathrm{b}}_{n-1} {\mathrm{b}}_n = 0$; 
we define $\Zrm_n(\Po,A) :=$ $\ker {\mathrm{b}}_n$,
$B_n(\Po,A) :=$ ${\mathrm{b}}_{n+1}(\Crm_{n+1}(\Po,A))$, and the 
{\em net homology group} $\Hrm_n(\Po,A) :=$ $\Zrm_n(\Po,A) / \Brm_n(\Po,A)$
(for $n=0$, we define $\Crm_0(\Po,A) := A$, $\Brm_0(\Po,A) := {\bf 0}$ and 
${\mathrm{b}}_0 := 0$, so that $\Hrm_0(\Po,A) = A$).
In the sequel, we will use
the notation ${\mathrm{b}} \equiv {\mathrm{b}}_n$, so that 
${\mathrm{b}}_{n-1} {\mathrm{b}}_n \equiv {\mathrm{b}}^2 = 0$.
Moreover, we will use the same notation for elements of
$\Zrm_n(\Po,A)$ and $\Hrm_n(\Po,A)$, identifying cycles
with the corresponding net homology classes.

When $A = \Zo$, $\Crm_n(\Po,\Zo)$ reduces to the Abelian group 
generated by the elements of $\widetilde \Si_n(\Po)$. 
Elements of $\Zrm_1(\Po,\Zo)$ satisfy the relation
\begin{equation}
\label{eq_z1}
\sum_i k_i \left( \partial_0 b_i - \partial_1 b_i \right) = 0
\ \ , \ \
\sum_i k_i b_i \in \Zrm_1(\Po,\Zo)
\ ,
\end{equation}
so that, $b \in \widetilde \Si_1(\Po) \cap \Zrm_1(\Po,\Zo)$ if, and only if,
$\partial_0b = \partial_1b$. 
We establish other useful relations. Given 
$a \in \Si_0(\Po)$  consider the degenerate $1$--simplex
$\sigma_0 a$  and note that 
$\sigma_0 a =$ ${\mathrm{b}} (\sigma_0\sigma_0 a)$.
So that,
\begin{equation}
\label{eq_s1_h1}
\sigma_0a \ = \ 0  \ \in \Hrm_1( \Po,\Zo) \ .
\end{equation}
If $b \in$ $\widetilde \Si_1(\Po)$ and $\overline b$ 
is defined by $\partial_0 \overline b :=$ $\partial_1 b$,
$\partial_1 \overline b :=$ $\partial_0 b$,
$|\overline b| :=$ $|b|$, then defining
$c :=$ $( |b| ; b , \sigma_0\partial_0b , \overline b )$ we obtain
\[
{\mathrm{b}} c = b - \sigma_0\partial_0b + \overline b
\ ,
\]
implying (by (\ref{eq_s1_h1}))
\begin{equation}
\label{eq_bbbar}
\overline b + b = 0 \ \in \Hrm_1(\Po,\Zo) \ .
\end{equation}
Now, for every $n \in \No$ we consider the map
\begin{equation}
\label{eq_derham}
\Crm^n ( \Po , A ) \times \Crm_n (\Po,A)  \to A
\end{equation}
obtained extending the evaluation $( v , x ) \mapsto v(x)$, 
$v \in \Crm^n(\Po,A)$, $x \in \widetilde \Si_n(\Po)$, by 
linearity. Some elementary computations show that
(\ref{eq_derham}) induces a bilinear map
\begin{equation}
\label{eq_derham2}
\Hrm^n (\Po,A) \times \Hrm_n(\Po,A) \to A
\ .
\end{equation}

\

We now analyze the connection between homotopy and net homology,
providing a version of a classical result (\cite[Thm.2.A.1]{Hat}).
To this end, we make some preliminary remarks on 
$\pi_1(\Po,a)$, $a \in \Si_0(\Po)$,
and the associated Abelianized group $\pi_1(\Po,a)_{ab}$.

We recall that $\pi_1(\Po,a)_{ab}$ is defined as the quotient
of $\pi_1(\Po,a)$ by the commutator subgroup, that is the normal subgroup 
generated by elements of the form $p_1 * p_2 * p_1^{-1} * p_2^{-1}$,
$p_1$, $p_2 \in \pi_1(\Po,a)$; by construction, $\pi_1(\Po,a)_{ab}$ is the 
universal Abelian group such that each group morphism $\pi_1(\Po,a) \to A$,
with $A$ Abelian, factorizes through a morphism $\pi_1(\Po,a)_{ab} \to A$.
For each $p \in \pi_1(\Po,a)$, we denote the corresponding class 
in $\pi_1(\Po,a)_{ab}$ by $[p]$, so that 
$[p*p'] = [p]*[p'] = [p']*[p] = [p'*p]$, $p,p' \in \Po (a)$.

\begin{oss}
\label{rem_cp_tr}
Let $c \in \widetilde \Si_2(\Po)$. Then the path
$p_c :=  \overline{\partial_1 c} * \partial_0c * \partial_2c$
is homotopic to the constant path.
\end{oss}

For each $a \in \Si_0(\Po)$, we denote the set of 
paths starting and ending in $a$ by $\Po (a)$.
\begin{oss}
\label{oss_p1}
Let $p' \in \Po(a')$, $a' \in \Si_0(\Po)$, $a' \neq a$; then there is 
a path $\gamma_{a,a'}$, starting in $a'$ and ending in $a$,
so that we can define 
$p := \gamma_{a,a'} * p * \overline \gamma_{a,a'} \in \Po(a)$.
\end{oss}

\begin{oss}
\label{rem_p1ab}
Let $p := b_n * \cdots * b_1 \in$ $\Po (a)$, $a \in \Si_0(\Po)$. 
Then for every index $i = 0 , \ldots , n$, we define the
path 
\[
p_i := b_i * \cdots * b_1 * b_n * \cdots * b_{i+1}
\in \Po(\partial_0 b_i)
\ .
\]
We say that $p_i$ is obtained from $p$ by a {\em shuffle}.
If $\gamma_{a,\partial_0b_i}$ is a path as in the previous remark, 
then we define
\[
\hat p_i :=
\gamma_{a,\partial_0b_i} * p_i * \overline \gamma_{a,\partial_0b_i} 
\in \Po(a) \ .
\]
In this way, we find (independently of the choice of 
$\gamma_{a,\partial_0b_i}$)
\[
\begin{array}{ll}
[\hat p_i ] & =
[ \gamma_{ a , \partial_0 b_i }  *
   b_i * \cdots * b_1  *
   b_n * \cdots * b_{i+1} *
  \overline \gamma_{ a , \partial_0 b_i } ] \\ & 
= 
[ \gamma_{ a , \partial_0 b_i } * b_i * \cdots * b_1 ] *
[ b_n * \cdots * b_{i+1} * \overline \gamma_{ a , \partial_0 b_i } ] \\ & 
= 
[ b_n * \cdots * b_{i+1} * \overline \gamma_{ a , \partial_0 b_i } ] *
[ \gamma_{ a , \partial_0 b_i } * b_i * \cdots * b_1 ]         \\ & 
= [p] 
\ .
\end{array}
\]
Of course, if $\partial_0 b_i = a$, then we can pick
$\gamma_{ a , \partial_0 b_i } = \sigma_0a$ and
$[p] = [p_i]$.
Note that if, for some $i$, $p_i$ is homotopic to
a constant path, then
$[ \hat p_i ] = 
 [ \gamma_{ a , \partial_0 b_i } * 
   \overline \gamma_{ a , \partial_0 b_i } ]$
and $[p]$ is the identity of $\pi_1(\Po,a)_{ab}$.
\end{oss}

\begin{lemma}
\label{lem_p1h1}
Let $a \in \Po$ and $p \in \Po (a)$ be a path of the form
\begin{equation}
\label{eq_deg_p}
p \ = \
  \ldots * \overline{b_m} * 
  \ldots * \overline{b_1} 
  \ldots * b_m * 
  \ldots * b_1 
\ \ , \ \
m = 2 , \ldots
\ .
\end{equation}
Then $[p] = [\sigma_0 a] \in \pi_1(\Po,a)_{ab}$.
\end{lemma}

\begin{proof}
$p$ is of the form
$
p  = p_2 * ( \overline{b_1} * p_1 * b_1 )$,
where $( \overline{b_1} * p_1 * b_1 ) , p_2 \in \Po (a)$.
Suppose that $p_1$ contains the $1$--simplex $b_m$ and $p_2$ $\overline b_m$.
Shuffle $\overline{b_1} * p_1 * b_1$ to give a path $p'_1$ ending with $b_m$; 
then, by Remark \ref{rem_p1ab}, 
$[\overline{b_1} * p_1 * b_1] = 
 [ \gamma_{a,\partial_0 b_m} * p'_1 * \overline \gamma_{a,\partial_0 b_m} ]$.
Shuffle $p_2$ to get a path $p'_2$ beginning with $\overline b_m$.
Then 
$[p_2] = 
 [\gamma_{a,\partial_0 b_m} * p'_2 * \overline \gamma_{a,\partial_0 b_m}]$.
It follows that 
\[
[p] = 
[ \gamma_{a,\partial_0 b_m} * p'_2 * p'_1 * \overline \gamma_{a,\partial_0 b_m} ]
\ . 
\]
Now, $p'_2 * p'_1$ contains $\overline b_m * b_m$ and $b_1 * \overline b_1$,
and these can be removed without changing the homotopy class; we then
have a path of the same type with fewer $1$--simplices.
Suppose on the other hand that there is no $b_m$ contained in $p_1$
such that $p_2$ contains $\overline b_m$. Then both $p_1$ and $p_2$ are
of the form (\ref{eq_deg_p}) and the result follows by induction.
\end{proof}

\begin{teo}
\label{teo_p1h1}
Let $( \Po , \leq )$ be a pathwise connected poset, and $a \in \Si_0(\Po)$.
Then there is an isomorphism $\Hrm_1(\Po,\Zo) \simeq$ $\pi_1(\Po,a)_{ab}$.
\end{teo}

\begin{proof}
{\bf Step 1}. Let $a \in \Si_0(\Po)$, and $p :=$ $b_n * \cdots * b_1$ a 
generic path in $\Po(a)$. With this notation, we define the group morphism
\begin{equation}
\label{eq_KH}
T : \Po (a) \to \Crm_1(\Po,\Zo)
\ \ , \ \
Tp := \sum_{i=0}^n b_i
\ .
\end{equation}
Since
\[
{\mathrm{b}} Tp = 
\sum_{i=0}^n  \left( \partial_0 b_i - \partial_1 b_i \right) =
a + \sum_{i=1}^{n-1} \left( \partial_0 b_i - \partial_1 b_{i+1} \right) - a
= 0
\ ,
\]
we conclude that $T$ actually takes values in $\Zrm_1(\Po,\Zo)$. 
{\bf Step 2}. We verify that (\ref{eq_KH}) factorizes through a map
\begin{equation}
\label{eq_KH1}
T : \pi_1 (\Po,a) \to \Hrm_1(\Po,\Zo)
\end{equation}
(note the slight abuse of the notation $T$). To this end, we consider
$c \in \widetilde \Si_2 (\Po)$, $p :=$ $b_n * \cdots * b_1 \in$ 
$\Po(a)$, and an elementary deformation 
of $p$ performed replacing a pair $b_i = \partial_0c$, $b_{i+1} = \partial_2 c$, 
$i \in \{ 0 , \ldots , n \}$, by $\partial_1c$. If we denote 
the deformed path by $p_c$, then
\[
Tp_c = 
\sum_{i=0}^n b_i - \sum_{k=0}^2 (-1)^k \partial_k c = 
Tp - {\mathrm{b}}_2c
\ .
\]
In the same way, an elementary deformation performed replacing 
$b_i = \partial_1 c$ by $\partial_0 c * \partial_2 c$ gives rise
to the operation
\[
Tp_c = 
\sum_{i=0}^n b_i + \sum_{k=0}^2 (-1)^k \partial_k c = 
Tp + {\mathrm{b}}_2c
\ .
\]
This implies that (\ref{eq_KH1}) is well-defined.
{\bf Step 3}.
We verify that (9) is surjective. To this end, note that a generic 
element of $\Zrm_1(\Po,\Zrm)$ can be written up to homology in the form 
\[
x = \sum_{j=1}^m b_j
\ ,
\]
no signs being needed since a $b_j$ can be replaced by $\overline b_j$ 
if necessary. We show that there is a bijection $f$ on $1,2,\dots,m$ and 
$0=m_0<m_1<\dots<m_\ell=m$ and paths 
$p_k=q_k*b_{f(m_k)}*b_{f(m_k-1)}*\dots*b_{f(m_{k-1}+1)}*\overline q_k\in \Po (a)$ 
for $k=1,2,\dots,\ell$. This suffices since setting $p:=p_\ell*\dots*p_1$, 
$Tp=x$ up to homology. We set $f(1)=1$ and since $x\in \Zrm_1(\Po,\Zrm)$ 
we may pick a $b_j$ with $\partial_0b_j=\partial_1b_1$ and set $f(2)=j$. 
We continue in this way, i.e.\ pick a $b_k$ with $\partial_0b_k=\partial_1b_j$ 
and set $f(3)=k$. After say $r$ steps this process terminates when two conditions 
are fulfilled, $\partial_1b_{f(r)}=\partial_0b_1$ and there is no $b_s$ with 
$\partial_0b_s=\partial_1b_{f(r)}$. We set $m_1:=r$ and pick a path $q_1$ with 
$\partial_0q=\partial_0b_1$ and $\partial_1q_1=a$. The sum over the remaining 
$1$--simplices is still in $\Zrm_1(\Po,\Zrm)$ and the argument may be repeated 
to reach the desired conclusion. Note that by Lemma \ref{lem_p1h1}, $Tp=Tp'$ 
implies $[p]=[p']$.
{\bf Step 4}. Let $p:=b_n*\cdots *b_1\in \ker T$, i.e.
$x := \sum_i b_i \in \Brm_1(\Po,\Zrm)$.
Then there are $2$--cycles $c_1,\dots , c_m$ such that 
\[
x = \sum_{j \in J} (\partial_0c_j - {\partial_1c}_j + \partial_2c_j) - 
    \sum_{k \in I} (\partial_0c_k - \partial_1c_k   + \partial_2c_k)
\ ,
\]
where $I \cup J = \{ 1 , \ldots , m \}$. Thus,
\[
\sum_i b_i+
\sum_j (\partial_1c_j + \overline{\partial_1c}_j) +
\sum_k (\partial_0c_k + \overline{\partial_0c_k} + 
        \partial_2c_k + \overline{\partial_2c_k})
=
\]
\[
\sum_j (\partial_0c_j  + \overline{\partial_1c}_j + \partial_2c_j ) +
\sum_k (\partial_0c'_k + \overline{\partial c}'_k + \partial_2c'_k) \ ,
\]
where $c'_k$ is the two simplex got by exchanging the vertices $0$ and $1$ of $c_k$. 
Adding on $\sum_j T(\gamma_j*\overline{\gamma}_j)+\sum_kT(\gamma_k*\overline{\gamma}_k)$,
where $\gamma_j$ is a path from $\partial_{01}c_j$ to $a$ and $\gamma_k$ a path from 
$\partial_{01}c_k$ to $a$, to each side, we get an element $y$ of $\Brm_1(\Po,\Zrm)$. 
From the form of the left hand side, we see that there is a path $p_1$ homotopic to $p$ 
with $Tp_1=y$ and, from the form of the right hand side, there is a path $p_2$ 
homotopic to $\sigma_0a$ with $Tp_2=y$. Hence $[p]=0$. 
\end{proof}

\begin{cor}
\label{cor_h1_p1}
Let $( \Po , \leq )$ be a pathwise connected poset. 
For every Abelian group $A$, there is an isomorphism 
$\Hrm^1(\Po,A) \simeq {\mathrm{Hom}} ( \Hrm_1(\Po,\Zo) , A )$.
\end{cor}

\begin{proof}
By Rem.\ref{oss_h1A}, we have an isomorphism 
$\Hrm^1(\Po,A) \simeq {\mathrm{Hom}} ( \pi_1(\Po,a) , A )$. 
Since $A$ is Abelian, each morphism 
$\phi \in {\mathrm{Hom}} ( \pi_1(\Po,a) , A )$
factorizes through an element of ${\mathrm{Hom}} ( \Hrm_1(\Po,\Zo) , A )$.
\end{proof}

Now, let $M$ be a Hausdorff space and let
$\Hrm_1(M,A)$, $\Hrm^1(M,A)$
denote respectively the singular homology and the singular cohomology
of $M$ with coefficients in $A$. We fix a poset $M_\prec$ 
providing a base of arcwise and simply connected open subsets of $M$. 
\begin{teo}
\label{teo_ho_co}
Let $M$ be a Hausdorff, locally arcwise and simply connected space. 
Then there is an isomorphism $\Hrm_1(M,\Zo) \simeq$ $\Hrm_1(M_\prec,\Zo)$. 
For every Abelian group $A$ there is an isomorphism 
$\Hrm^1(M_\prec,A) \simeq$ ${\mathrm{Hom}} ( \Hrm_1(M,\Zo) , A )$, 
and this implies the isomorphism
\begin{equation}
\label{eq_h1h1}
\Hrm^1(M_\prec,A) \simeq \Hrm^1(M,A)
\ .
\end{equation}
\end{teo}

\begin{proof}
Let $a \in \Si_0(M_\prec)$. 
By \cite[Thm.2.18]{Ruz05}, we have isomorphisms 
$\Hrm_1(M,\Zo) \simeq$ $\pi_1(M)_{ab} \simeq$ $\pi_1(M_\prec,a)_{ab}$
and, by the classical Hurewicz theorem there is an isomorphism
$\Hrm_1(M,\Zo) \simeq \pi_1(M)_{ab}$.
The isomorphism (\ref{eq_h1h1}) follows since, by the universal 
coefficient theorem, $\Hrm^1(M,A) \simeq {\mathrm{Hom}}( \Hrm_1(M,\Zo) , A )$ 
(see \cite[\S 3.1]{Hat}). % pag.198
\end{proof}

When $M$ is a manifold it makes sense to consider the
de Rham cohomology $\Hrm_{dR}^1(M,\Ro)$, and we have
\begin{cor}
Let $M$ be a connected manifold. Then there is an isomorphism
$\Hrm^1(M_\prec,\Ro) \simeq \Hrm_{dR}^1(M,\Ro)$.
\end{cor}

%*************************************************************

\section{Homotopy and net cohomology of net bundles.}

\subsection{Net structures.}
A {\em quasinet bundle} with {\em base} $\Po$ is given by a $4$-ple $\Xcal :=$ 
$(X,\pi,J,\Po)$, where $X$ is a set called the {\em total space},
$\pi : X \to \Po$ is a surjective map with {\em fibres} $X_a :=$ $\pi^{-1}(a)$,
$a \in \Si_0(\Po)$, and $J$ is a family of injective maps
$J_b : X_{\partial_1b} \to X_{\partial_0b}$,
$b \in \Si_1(\Po)$,
called {\em the net structure} of $\Xcal$, satisfying the {\em cocycle relations}
$J_{\partial_0c} J_{\partial_2c} =$ $J_{\partial_1c}$, $c \in \Si_2(\Po)$,
and $J_{\sigma_0a} =$ $id_{X_a}$, $a \in \Si_0(\Po)$. 
When each $J_b$ is a bijective map, we say that $\Xcal$ 
is a {\em net bundle}. \smallskip

% In order to economize on notation, 
% if $b := (\partial_1b,\partial_0b) \in \Si_1(\Po)$ and 
% $\overline b := (\partial_0b,\partial_1b) \in \Si_1(\hat \Po)$, 
% then we write
% %
% \begin{equation}
% \label{eq_J}
% J_b := J_{(\partial_0b,\partial_1b)}
% \ \ , \ \
% J_{\overline{b}} := J_b^{-1} 
% \ .
% \end{equation}
% %

The fibres of a net bundle $\Xcal$ are all isomorphic
(see \cite{RRV07}); a distinguished 
fibre $X_a$ of $\Xcal$ will be called {\em the 
standard fibre of} $\Xcal$ and emphasized with the notation $F$. 

Analogous definitions apply when considering categories with more structure
than the category of sets; 
in that case, the net structure shall be a family of monomorphisms (or isomorphisms)
in the appropriate category. Of particular interest for our purposes will be
{\em Hilbert (quasi)net bundles, group (quasi)net bundles, {\it C*}-algebra 
(quasi)net bundles and (quasi)net bundles of topological spaces}.

The first class of quasinet bundles will be studied in \S \ref{sec_nkt}. 
Quasinet bundles of {\it C*}-algebras motivated  our work. They arise as follows.

\begin{ex}
In the algebraic approach to quantum field theory (see \cite{DR90}, for example), 
the set of quantum observables 
is presented as a {\em net of {\it C*}-algebras}; by this we mean a map 
$\mathcal{A}$, assigning a {\it C*}-algebra $\mathcal{A}(\dc)$ to elements $\dc$ 
of a distinguished poset $\Po$ of contractible open subsets of spacetime, and 
interpreted as the algebra of quantum observables localized in the open set 
$\dc \in \Po$. The map $\mathcal{A}$ preserves the order, in the sense that 
$\mathcal{A}(\dc_1) \subseteq$ $\mathcal{A}(\dc_2)$ 
whenever $\dc_1 \subseteq \dc_2$.
This structure may be replaced by a quasinet bundle of {\it C*}-algebras
defined as follows: we consider the set
\[
\hat{\mathcal A} 
:= 
\left\{ 
( \dc , \mathcal{A}(\dc) ) 
\ , \ \dc \in \Po
\right\}
\]
endowed with the projection $\pi : \hat{\mathcal A} \to \Po$ 
onto the first component and net structure $J_b$, $b \in \Si_1(\Po)$, 
defined by the inclusion 
$\mathcal{A}(\partial_1b) \subseteq \mathcal{A}(\partial_0b)$.
\end{ex}

Now let $\hat \Xcal :=$ $(\hat X, \hat \pi, \hat J,\Po )$ be a quasinet bundle.
A map $T : X \to X$ is said to be a {\em morphism} if $\hat \pi T = \pi$ (this
implies that $T$ restricts to maps $T_a : X_a \to \hat X_a$, $a \in \Si_0(\Po)$) 
and
\[
\hat J_b T_{\partial_1b} = T_{\partial_0b} J_b
\ \ , \ \
b \in \Si_1(\Po)
\ .
\]
In this case, we use the notation $T \in ( \Xcal , \hat \Xcal )$; if each $T_a$, 
$a \in \Si_0(\Po)$, is a bijective map, then we say that
$T$ is an {\em isomorphism}. In this way, the set of
quasinet bundles becomes a category. The full subcategory of net bundles with
fibre $F$ will be denoted by
$\mB ( \Po , F )$.

Every net bundle $\Xcal :=$ $( X, \pi, J,\Po )$ defines a 
net bundle  
\begin{equation}
\label{def_Xc}
\Xcal_\circ := ( X_\circ , \pi_\circ , J^\circ , \Po^\circ )
\ ,
\end{equation}
where   $\Po^\circ$ is the dual poset
(see Appendix \ref{sec_ss_poset}), 
$X_\circ := X$, $\pi_\circ := \pi$, and the net structure
$J^\circ$ is defined by $J^\circ_b := J^{-1}_{(\partial_1b,\partial_0b)}$, 
$b \in \Si_1(\Po^\circ)$. 
This provides an isomorphism $\mB ( \Po , F ) \simeq \mB ( \Po^\circ , F )$,
which has been described in cohomological terms in \cite[\S 7]{RRV07}.

Now, let $G$ denote the group of automorphisms of $F$. 
By \cite[Prop.3.8]{RR06}, \cite[Thm.6.7]{RRV07}, 
for each $a \in \Si_0(\Po)$ there are isomorphisms
\begin{equation}
\label{eq_b_c}
\dot{\mB} ( \Po , F ) 
\simeq 
\Hrm^1 ( \Po , G ) 
\simeq 
\dot{\mathrm{H}}{\mathrm{om}} (\pi_1(\Po,a),G)
\ ,
\end{equation}
where $\dot{\mB} ( \Po , F ) $ is the set of isomorphism classes of net bundles 
with fibre $F$ and $\dot{\mathrm{H}}{\mathrm{om}} (\pi_1(\Po,a),G)$ denotes
the set of equivalence classes of morphisms from the homotopy group 
$\pi_1(\Po,a)$ to $G$, two morphisms being equivalent if they
differ by an inner automorphism of $G$.

The notion of morphism can be generalized to allow changes 
of the base $\Po$; to this end, we introduce the notion of {\em pullback}.

Let $\eta : \Po' \to \Po$ be a morphism of posets. Given a quasinet bundle 
$\Xcal$ over $\Po$, we define a quasinet bundle $\Xcal^\eta$, by considering 
the set
\[
X^\eta := X \times_\Po \Po' :=
\left\{ 
(x,a') \in X \times \Si_0(\Po') :  \pi(x) = \eta (a')
\right\}
\]
endowed with the obvious projection $\pi^\eta : X^\eta \to \Po'$ and the net structure
$J^\eta_{b'} (x,\partial_1b') :=$ $( J_{\eta(b')}x , \partial_0b' )$, 
$x \in X_{\partial_1 \eta (b')}$.
By restricting the notion of pullback to the subcategory of net bundles with 
fibre $F$, we obtain a functor
\[
\eta^* : \mB ( \Po , F ) \to \mB ( \Po' , F )
\ \ , \ \
\eta^* \Xcal := \Xcal^\eta
\ .
\]
We are now able to define the notion of a {\em generalized morphism} from a quasinet
bundle $\hat \Xcal :=$ $( \hat X , \hat  \pi , \hat J , \hat \Po )$ to a quasinet bundle
$\Xcal :=$ $( X , \pi , J , \Po )$ as a pair $(\eta , T)$, where
$\eta : \hat \Po \to \Po$ is a morphism of posets, and $T \in ( \Xcal^\eta , \hat \Xcal )$ 
is a morphism in the usual sense.

We also recall that a {\em local section} of a quasinet bundle $( X , \pi , J , \Po )$
is given by a map $\sigma$ from a subposet $V \subseteq \Po$ into $X$, such that $\pi \sigma =$
$id_V$ and $J_b \sigma (\partial_1b) =$ $\sigma(\partial_0b)$, $b \in \Si_1(V)$.
If $V = \Po$, then $\sigma$ is said to be {\em global}. We denote 
the set of local sections of $\Xcal$ over $V$ by
$S ( V ; \Xcal )$
(for $V = \Po$, we use the analogous notation for the set of global sections).

\subsection{Poset net bundles.}
\label{pullback}

In the present section, we consider the notion of poset net bundle. The interesting 
applications are when the fibre is the poset arising from the topology of a space. 
In particular, in the following we relate the homotopy of a poset net bundle
to that of the underlying poset.

A {\em poset net bundle} is a net bundle $( X , \eta , J , \Po )$ where each fibre $X_a$,
$a \in \Si_0(\Po)$, is endowed with an ordering $\leq_a$ and each map $J_b$,
$b \in \Si_1(\Po)$, is an isomorphism of posets.

A standard argument for net bundles allows one to conclude that each poset net bundle 
admits a fixed poset $( F , \leq )$ as standard fibre, endowed with poset isomorphisms 
$V_a : F \to X_a$, $a \in \Si_0 (\Po)$.

\begin{lemma}
\label{lem_ord}
Let $(X,\eta,J,\Po)$ be a poset net bundle. Then there is a canonical ordering 
$\prec$ on the total space $X$ making $\eta$ a poset epimorphism. In this way, 
morphisms of poset net bundles give rise to poset morphisms of the underlying 
total spaces.
\end{lemma}
\begin{proof}
Define $x_1 \prec x_0$ $\Leftrightarrow$ $a_1 := \eta(x_1) \leq \eta(x_0) =: a_0$ 
and $J_{(a_0,a_1)} (x_1) \leq_{a_0} x_0$. 
\end{proof}
The previous lemma illustrates the main advantage of considering poset net bundles: 
in fact, the usual net bundles $( X , p , J , \Po )$ correspond to partial orderings 
where each fibre $X_a$ has the ``discrete'' order relation 
$x \leq x'$ $\Leftrightarrow$ $x = x'$,
$x,x' \in X_a$ (see \cite[\S 4.1]{RRV07}), whereas a poset net bundle may be 
endowed with a more interesting partial ordering.

By the previous lemma it makes sense to consider the maps
\begin{equation}
\label{def_pbC}
\eta^* : \Hrm^1 ( \Po , G ) \rightarrow \Hrm^1 ( X , G ) 
\ \ ,
\end{equation}
\begin{equation}
\label{def_ppi1}
\eta_* : \pi_1 (X,x) \rightarrow \pi_1 (\Po,a)
\ \ , \ \
x \in X_a 
\ ,
\end{equation}
for every poset net bundle $( X , \eta , J , \Po )$ and 
group $G$.

\begin{teo}
\label{teo_exp1}
Let $( X , \eta , J , \Po )$ be a poset net bundle with pathwise connected fibre
$F$. Then for each $a \in \Si_0(\Po)$ there is a morphism
$j_a : \pi_1(X_a,x) \to \pi_1(X,x)$ with $\pi_1(X_a,x) \simeq \pi_1(F,a')$, 
$x \in X_a$, $a' \in \Si_0(F)$,
and an exact sequence
\begin{equation}
\label{eq_ex_p}
\pi_1 (X_a,x)
\stackrel{j_a}{\to}
\pi_1 (X,x)
\stackrel{\eta_*}{\to}
\pi_1 (\Po,a)
\to
{\bf 0}
\ \ .
\end{equation}
In particular, if $(F,\leq)$ is simply connected, then $\eta_*$ is an isomorphism
\end{teo}

\begin{proof}
The morphism $j_a$ is induced by the order-preserving inclusion
$j_a :( X_a , \leq_a  ) \to ( X , \prec  )$.
Clearly,
\[
\eta_* j_a [p] \ = \ [\sigma_0 a]
\ \ , \ \
[p] \in \pi_1(X_a,x)
\ ,
\]
where $\sigma_0 a$, $a \in \Si_0 (\Po)$, denotes, as usual, the degenerate 
$1$-simplex. In other words, $j_a (\pi_1 (X_a,x)) \subseteq$ $\ker \eta_*$.
Let $p :=$ $b_n * \cdots * b_1$ be a closed path in $\Si_1(\Po)$.
We fix $|v_1| \in X_{|b_1|}$ and define recursively
\[
\left\{
\begin{array}{ll}
|v_{k+1}|  :=
J_{|b_{k+1}|,\partial_0 b_k} J_{|b_k|,\partial_0 b_k}^{-1}
|v_k|
\\
\partial_i v_{k+1} :=
J_{|b_{k+1}|,\partial_i b_{k+1} }^{-1}  |v_{k+1}|
\ , \
i = 0,1
\\
v_k :=
(  
|v_k| ;
\partial_0 v_k ,
\partial_1 v_k
)
\end{array}
\right.
\]
(recalling that $\partial_0 b_k = \partial_1 b_{k+1}$, $k = 1 , \ldots , n-1$). 
In this way, we obtain a path $\tilde p_0 :=$ 
$\widetilde b_n * \cdots * \widetilde b_1$ in $\Si_1(X)$
such that  $\eta_* \tilde p_0 =$ $p$. Since $\tilde p_0$ may be not closed,
we consider a path $\tilde p_1$ in $\Si_1( X_{\partial_1 b_1} \simeq F )$ 
starting in $\partial_0 v_n$ and ending in $\partial_1 v_1$, and 
define $\tilde p := \tilde p_0 * \tilde p_1$. By construction, 
$\eta_* (\tilde p)$ coincides with the closed path
\[
b_n * \ldots * b_1 *
\sigma_0 \partial_1 b_1 * \ldots * \sigma_0 \partial_1 b_1
\]
where $\sigma_0 \partial_1 b_1$ is the degenerate 1-simplex in $\Si_1 (\Po)$. 
The closed path constructed above is clearly homotopic to $p$.
For the converse, let $\tilde p :=$ $v_n * \cdots * v_1$ be a closed
path in $X$ with $\tilde p \in \ker \eta_*$. We prove that $\tilde p$ is homotopic
to a path in $X_a$ for some $a \in \Si_0(\Po)$. To this end, we consider the
path in $\Po$
\[
\eta_*(\tilde p)
:=
b_n * b_{n-1} * \cdots * b_1 
\ , \
b_k := \eta_* v_k
\ , \
k = n , \ldots , 1
\ .
\]
Since $\tilde p \in \ker \eta_*$, we may find $2$-simplices $c_1 , \ldots , c_m \in$
$\Si_2(\Po)$ providing elementary deformations of $\eta_* (\tilde p)$ to the
constant path
\[
\sigma_0 a
\ \ , \ \
a := \partial_1 b_n
\ , \
\sigma_0 a := ( a; a,a )
\ .
\]
In particular we may have, for example, an elementary deformation
contracting $b_2 * b_1$ to $\partial_1 c_1$.
Now, $v_2 * v_1$ is a path in $\eta^{-1}(V_{|c_1|})$;
using the definition of local chart (\cite[\S 4.3]{RRV07}) and
Lemma \ref{lem_ord}, we conclude that there is an isomorphism
of posets
\[
\eta^{-1}(V_{|c_1|}) \simeq V_{|c_1|} \times F
\ .
\]
The image of $v_2 * v_1$ under the above isomorphism 
can be written, according to \S \ref{sec_hom_prod}, as a pair $(p^\alpha,p^\mu)$, 
where $p^\alpha = \partial_1 c_1$ is a path in $V_{|c_1|}$ and $p^\mu$ is a
path in $F$. Applying (\ref{eq_pss}), we conclude that $v_2 * v_1$
is homotopic to a path $v'_m * \ldots * v'_1$ such that
\[
\eta_* (v'_m * \ldots * v'_1) 
= 
\partial_1 c_1 * \sigma_0 a * \ldots * \sigma_0 a
\ .
\]
Thus, $\tilde p$ is homotopic to a path of the type
\[
\tilde p_1 
\ := \ 
v_n * v_{n-1} * \cdots * v_3 * v'_m * \ldots * v'_1
\ ,
\]
Note that if there is a $c_2 \in \widetilde \Si_2(\Po)$ contracting 
$b_3 * \partial_1 c_1$, then $v_3 * v'_m * \ldots * v'_1$ 
is a path in $\eta^{-1} (V_{|c_2|}) \simeq V_{|c_2|} \times F$, so that 
we can again apply  the above argument.
Moreover, the same procedure applies to ampliations of $\eta_*(\tilde p)$.
Iterating the above operations for each deformation of $\eta_*(\tilde p)$,
we conclude that $\tilde p$ is homotopic to a path $\tilde p_m$ such that 
$\eta_*(\tilde p_m) = \sigma_0 a * \ldots * \sigma_0 a$, 
i.e. $\tilde p_m$ is a path in $X_a$. 
This proves $\ker \eta_* = j_a \left( \pi_1 (X_a,x) \right)$.
\end{proof}

\begin{cor}
\label{cor_CV}
Let $( X , \eta , J , \Po )$ be a poset net bundle with simply 
connected fibre $(F,\leq)$. Then the maps
(\ref{def_pbC}), (\ref{def_ppi1}) define group isomorphisms.
\end{cor}
\begin{proof}
It suffices to note that $\pi_1(F,a') = {\bf 0}$
and apply (\ref{eq_b_c}), (\ref{eq_ex_p}).
\end{proof}

The previous theorem has a well-known counterpart, 
involving topological spaces; in this case, (\ref{eq_ex_p})
is a long exact sequence involving higher homotopy groups.
One might hope to generalize Thm.\ref{teo_exp1} and get 
a long exact sequence by introducing higher homotopy groups for posets.

\

Now let $\Xcal :=$ $( X , p , J , \Po )$ be a net bundle of topological spaces, 
so that each $J_b$, $b \in \Si_1(\Po)$, is a homeomorphism.
We fix a poset $X_{\prec,a}$ of open subsets $U \subset X_a$, $U \neq \emptyset$,
ordered under inclusion and forming a base for $X_a$.
We require that each $J_b$ is a poset
isomorphism from $X_{\prec,\partial_1b}$ to $X_{\prec,\partial_0b}$. 
Then defining 
$X_{\prec} :=$ $\dot{\cup}_a X_{\prec,a}$ and letting 
$p_{\prec} : X_{\prec} \to \Po$ be the obvious projection, 
we conclude that
\begin{equation}
\label{eq_Xprec}
\Xcal_{\prec} := ( X_{\prec} , p_{\prec} , J , \Po  )
\end{equation}
is a poset net bundle, called the {\em poset net bundle associated with} 
$\Xcal$.
By \cite[Thm.2.8]{Ruz05}, the topological homotopy group $\pi_1(X_a)$
is isomorphic to $\pi_1 ( X_{\prec,a} , U )$, 
$a \in \Si_0(\Po)$, $U \in \Si_0(X_{\prec,a})$;
thus, by Thm.\ref{teo_exp1} we conclude the following:
\begin{cor}
\label{cor_top}
Let $\Xcal$ be a net bundle of topological spaces with fibre a Hausdorff,
locally arcwise and simply connected space $M$. For each $a \in \Si_0(\Po)$, 
there is an exact sequence
\begin{equation}
\label{eq_ex_t}
\pi_1 (M)
\stackrel{j_a}{\to}
\pi_1 (X_{\prec},U)
\stackrel{\eta_*}{\to}
\pi_1 (\Po,a)
\to
{\bf 0}
\ \ .
\end{equation}
\end{cor}

\begin{ex}
For the notion of principal net bundle, we refer the reader to \cite{RRV07}.
Let $\Pcal :=$ $( P , \pi , J , \Po , R )$ be a principal net bundle with
a connected Lie group $G$ as standard fibre. Then (provided the net
structure is defined by means of isomorphisms of Lie groups), $\Pcal$ is
also endowed with the structure of a topological
net bundle and the previous corollary applies.
\end{ex}

The global space of a net bundle of topological spaces $\Xcal:=(X,p,J,\Po)$ can be 
topologized in the following way: pick for each $a\in\Si_0(\Po)$ a base $\Ucal_a$ 
for the topology of $X_a$. If $U\in\Ucal_a$, we define the {\it cylinder with base} 
$U$ to be 
\begin{equation}
\label{eq_cyl}
T_{a,U}:=\{x\in J_{(\dc,a)}(U), \dc \in V_a\} \ .
\end{equation}
Clearly $p(T_{a,U})=V_a$. The family $\{T_{a,U}\}$ is a base for a topology on $X$ 
and we denote the associated topological space by $\tau X$. This topology is independent 
of the choice of bases. Note that there is a continuous bijection 
\[
T_{a,U}\to V_a\times U
\ \ , \ \ 
x \mapsto ( p(x) , J_{(a,p(x)}x ) 
\ .
\]
The projection $p$ is therefore continuous as a map from $\tau X$ to $\tau K$.

\begin{lemma} 
\label{lem_net_tau}
Let $\Xcal := (X,p,J,\Po)$, $\hat\Xcal := (\hat X,\hat p,\hat J,\Po)$ be net bundles 
of topological spaces and $f \in (\Xcal,\hat \Xcal)$ a morphism. Then 
$f:\tau X\to\tau\hat X$ is continuous.
\end{lemma}

\begin{proof} 
It suffices to show that if $\hat V$ is open in $\hat X_a$ then 
$f^{-1}(T_{a,\hat V})=T_{a,f^{-1}(\hat V)}$. If 
$x\in J_{(\dc,a)}(f^{-1}(\hat V))$ then 
$f(x)\in fJ_{(\dc,a)}(f^{-1}(\hat V))\subset \hat J_{(\dc,a)}(\hat V)$. 
Hence $T_{a,f^{-1}(\hat V)}\subset f^{-1}(T_{a,\hat V})$. 
If $f(x)\in \hat J_{(\dc,a)}(\hat V)$ then 
$x\in f^{-1}\hat J_{\dc,a}(\hat V)=   J_{\dc,a} f^{-1} (\hat V)$ 
so $f^{-1}(T_{a,\hat V})\subset T_{a,f^{-1}(\hat V)}$.
\end{proof}

When $\Xcal$ is trivial, the previous Lemma implies that
there is a homeomorphism
\begin{equation}
\label{eq_tt}
{\tau X} \simeq {\tau \Po} \times M
\ ,
\end{equation}
where $\tau \Po$ is the space $\Po$ with the Alexandroff
topology (see \cite[\S 2.3]{RRV07}). 
The previous elementary remark gives an idea of the local behaviour 
of ${\tau X}$ for a generic $\Xcal$, and implies that the bundle 
$p : \tau X \to \tau \Po$ is locally trivial when $\Xcal$ is locally
trivial.

If $s \in S ( V_a ; \Xcal)$ and $s(a) \in U$, then 
$s(\dc) \in T_{a,U}$ for each $\dc \in V_a$ and 
$s^{-1}(T_{a,U}) = V_a$; we conclude that each 
$s \in S (V_a ; \Xcal)$ defines a local section 
$s : V_a \to {\tau X}$.

We now consider the associated net bundle of topological spaces 
$\Xcal_{\circ}$, see  (\ref{def_Xc}), and the associated poset 
$X_{\circ,\prec}$, see (\ref{eq_Xprec}) and 
Lemma \ref{lem_ord}.
\begin{teo}
\label{teo_net_top}
Let $\Xcal :=$ $( X , p , J , \Po )$ be a net bundle of topological spaces 
with standard fibre a space $M$. Then the map $p : {\tau X} \to {\tau \Po}$ 
defines a fibre bundle with fibre $M$. Moreover, there is an isomorphism
\begin{equation}
\label{eq_is_t_p}
\pi_1( X_{\circ,\prec} , U ) \simeq \pi_1({\tau X})
\ \ , \ \
U \in \Si_0(X_{\circ,\prec})
\ ,
\end{equation}
and for every group $G$ there are isomorphisms
\begin{equation}
\label{eq_t_h1}
\Hrm^1( X_{\circ,\prec} , G)
\simeq
{\mathrm{\dot{H}om}} ( \pi_1({\tau X}) , G)
\ .
\end{equation}
\end{teo}
\begin{proof}
In order to prove the theorem it just remains to check the details 
in (\ref{eq_is_t_p}).
To this end, let us consider the poset $X_{cyl}$ with elements the
sets (\ref{eq_cyl}), ordered under inclusion. 
Since $X_{cyl}$ is a base for $\tau X$,
we have an isomorphism $\pi_1(X_{cyl},T) \simeq \pi_1(\tau X)$,
$T \in \Si_0(X_{cyl})$.
Thus, in order to prove (\ref{eq_is_t_p}), it suffices to construct
a poset isomorphism from $X_{\circ,\prec}$ to $X_{cyl}$.
By the definition of $\Xcal_\circ$, $X_{\circ,\prec}$
coincides with $X_\prec$ as a set, 
thus elements of $X_{\circ,\prec}$ are given by open sets 
$U \in X_{\prec,a}$, $a \in \Si_0(\Po^\circ)$. 
The order relation for $X_{\circ,\prec}$ is given by 
$U \prec U'$ $\Leftrightarrow$ $a \geq a'$ and
$J^\circ_{(a',a)}(U) \subseteq U'$, i.e.
$U \subseteq J_{(a,a')}(U')$ (see Lemma \ref{lem_ord}).
We consider the bijective map
\begin{equation}
\label{eq_Td}
X_\prec \to X_{cyl}
\ \ , \ \
U \in X_{\prec,a} \mapsto T_{a,U}
\ .
\end{equation}
By construction, $U \prec U'$ implies $V_a \subseteq V_{a'}$
and $U \subseteq J_{(a,a')}(U')$, thus 
$T_{a,U} \subseteq T_{a',U'}$, and the theorem is proved.
\end{proof}

Let $M$ be a space and $P \to M$ a locally trivial bundle. 
We say that $P$ is \emph{locally constant} whenever it has a set 
of locally constant transition maps (that is, the transition maps
are constant on the connected components of their domains).
\begin{cor}
\label{cor_h1_net_top}
Let $G$ be a group. For every principal net $G$-bundle $\Pcal :=$
$(  P,p,J,R,\Po  )$, the fibre bundle $p : \tau P \to {\tau \Po}$
is locally constant.
\end{cor}

\begin{proof}
For every $a \in \Si_0(\Po)$, $V_a$ is an open neighbourhood of
$a$ trivializing $\tau P$, so there are local charts
$\theta_a : V_a \times G \to$ $p^{-1}(V_a)$, $a \in \Si_0(\Po)$,
giving rise to locally constant transition maps 
(see \cite[Lemma 5.6]{RRV07}).
\end{proof}

\section{K-theory of a poset.}
\label{sec_nkt}
\subsection{Net bundles of Banach spaces.}

We now specialize our discussion to the case of quasinet bundles of Banach spaces.
These objects have already been studied in a combinatorial setting (\cite{L}).
In the present paper, we emphasize the topological viewpoint.

Banach quasinet bundles are quasinet bundles having Banach spaces as fibres, 
and net structure given by bounded, injective, linear operators.
If $\Ecal := ( E , \pi , J , \Po )$, 
$\hat \Ecal := ( \hat E , \hat \pi , \hat J , \hat \Po )$ 
are Banach quasinet bundles, then we denote the set of bounded morphisms 
from $\Ecal$ into $\hat \Ecal$ by $( \Ecal , \hat \Ecal )$.
If $T \in$ $( \Ecal , \hat \Ecal )$, then each $T_a : E_a \to \hat E_a$,
$a \in \Si_0(\Po)$, is a bounded linear operator satisfying the relations
\begin{equation}
\label{eq_TbJ}
T_{\partial_0b} J_b = \hat J_b T_{\partial_1b} \ \ , \ \ b \in \Si_1(\Po) \ .
\end{equation}
If $\Ecal, \hat \Ecal$ are Banach net bundles, then
($\Po$ being pathwise connected)
(\ref{eq_TbJ}) implies that $\left\| T_a \right\|$ does not depend
on the choice of $a \in \Si_0(\Po)$.

The Banach quasinet bundle $\Ecal$ is said to be {\em finite dimensional} 
if the dimension of the fibres $E_a$, $a \in \Si_0(\Po)$, has an upper 
bound $d \in \No$. Since each $J_b$, $b \in \Si_1(\Po)$, between
vector spaces with the same finite dimension is also surjective,
we have

\begin{lemma}
Let $\Ecal$ be a finite-dimensional Banach quasinet bundle. 
Then $\Ecal$ is a Banach net bundle if and only if the rank 
function $d(a) :=$ ${\mathrm{dim}} (E_a)$, $a \in \Si_0(\Po)$, 
is constant.
\end{lemma}

A Banach quasinet bundle $\Ecal$ is said to be a {\em Hilbert quasinet bundle} 
if each fibre $E_a$, $a \in \Si_0(\Po)$, is a Hilbert space, 
and each $J_b$, $b \in \Si_1(\Po)$, preserves the scalar product.

Some remarks on the classification of Banach net bundles follow.
If we want to apply (\ref{eq_b_c}) to Banach net bundles,
then we have to pick the group $G$ to take account of the
structure of interest. For example, we may consider
the group of invertible operators of a Banach space, or the unitary group
$\Ucal$ of a Hilbert space $\Hcal$ if we are interested in isomorphisms 
preserving the Hermitian structure.
In particular, in the finite-dimensional case, 
$G$ is the complex linear group $\mathbb{GL}(d)$,
$d \in \No$, or the unitary group $\Uo(d)$.

Let us denote the set of isomorphism classes of Hilbert net bundles 
with fibre $\Hcal$ by $\dot{H}_{net}(\Po,\Hcal)$, and the \v Cech 
cohomology of the dual poset $\Po^\circ$ by $\Hrm^1_c(\Po^\circ,\Ucal)$. 
Applying \cite[Thm.8.1]{RRV07}, we find
\begin{prop}
\label{prop_vs_sc}
For each poset $\Po$ and $a \in \Si_0(\Po)$, there are 
isomorphisms
\[
\dot{H}_{net}(\Po,\Hcal) \simeq
\Hrm^1(\Po,\Ucal) \simeq
\Hrm^1_c(\Po^\circ,\Ucal) \simeq
\dot{\mathrm{H}}{\mathrm{om}} (\pi_1(\Po,a),\Ucal)
\ .
\]
\end{prop}

Unlike ordinary vector bundles, a finite-dimensional Banach net bundle
cannot in general be endowed with an Hermitian structure, in fact the 
equivalence relation induced by inner automorphisms leaves the 
determinant map
$\det \chi (p)$,
$\chi \in {\mathrm{Hom}} ( \pi_1(\Po,a) , \mathbb{GL}(d) )$,
$p \in \pi_1(\Po,a)$,
invariant. Thus, if $\chi$ does not take values in $\Uo(d)$, the same 
is true of every $\chi' \in {\mathrm{Hom}}(\pi_1(\Po,a),\mathbb{GL}(d))$
equivalent to $\chi$.

\

We denote the category of Hilbert quasinet bundles over $\Po$ by 
$qH_{net}(\Po)$, and the full subcategory of Hilbert net bundles by 
$H_{net}(\Po)$. For finite dimensional Hilbert quasinet bundles, 
we use the analogous notations $qV_{net}(\Po)$ and $V_{net}(\Po)$.

Some algebraic structures are naturally defined on $qH_{net}(\Po)$:
(1) The direct sum $\Ecal \oplus \hat \Ecal$;
(2) The adjoint $* : ( \Ecal , \hat \Ecal ) \to ( \hat \Ecal , \Ecal )$;
(3) The tensor product $\Ecal \otimes \hat \Ecal$;
(4) The symmetry $\theta \in$ 
$(\Ecal \otimes \hat \Ecal , \hat \Ecal \otimes \Ecal)$,
$\theta_a (v \otimes \hat v) :=$ $\hat v \otimes v$, 
$v \in E_a$, $\hat v \in \hat E_a$;
(5) The conjugate net bundle $\overline \Ecal :=$ 
$( \overline E , \pi , \overline J , \Po )$;
(6) The antisymmetric tensor powers $\lambda^r \Ecal$, $r \in \No$.

There is a Banach quasinet bundle of morphisms
$B(\Ecal,\hat \Ecal)$ associated with Hilbert quasinet bundles
$\Ecal$, $\hat \Ecal$. It is defined as follows:
for each $a \in \Po$, consider the vector space $( E_a , \hat E_a )$ of
linear operators from $E_a$ into $\hat E_a$ and the disjoint union $B(E,E') :=$
$\dot \cup_a ( E_a , \hat E_a )$. Then define the net structure $I_b (t) :=$
$\hat J_b  t  J_b^{-1}$, $t \in (E_{\partial_1b},\hat E_{\partial_1b})$,
$b \in \widetilde \Si_1(\Po)$.
In the finite-dimensional case, we clearly have an 
isomorphism of Banach quasinet bundles
\begin{equation}
\label{eq_iso_bm}
B(\Ecal,\hat \Ecal) \to \overline \Ecal \otimes \hat \Ecal
\ .
\end{equation}
The above isomorphism induces an Hermitian structure on $B(\Ecal,\hat \Ecal)$,
making it into a Hilbert quasinet bundle.

\begin{lemma}
\label{lem_sbm}
Every morphism $T \in ( \Ecal , \hat \Ecal )$ defines 
a section of $B(\Ecal,\hat \Ecal)$.
\end{lemma}

\begin{proof}
If $T \in ( \Ecal , \hat \Ecal )$, then we may regard $T$ as a map
$T : \Si_0(\Po) \to B(\Ecal,\hat \Ecal)$, $T(a) :=$ $T_a$. By the 
definition of morphism, we find that $T$ satisfies the compatibility 
property w.r.t. the net structure of $B(\Ecal,\hat \Ecal)$, 
i.e. $I_b \circ T (\partial_1b) =$ $T(\partial_0b)$.
\end{proof}

The existence of the adjoint and norm on $qH_{net}(\Po)$ imply
that the space $( \Ecal , \Ecal )$ is a unital {\it C*}-algebra. 
We can now prove the following:
\begin{lemma}
The category $qH_{net}(\Po)$ has subobjects.
\end{lemma}

\begin{proof}
Consider $\Ecal \in qH_{net}(\Po)$, $\Ecal :=$ $( E, \pi , J , \Po )$ 
and a projection $P \in ( \Ecal , \Ecal )$. We define 
$P \Ecal := ( PE , \pi' , J' , \Po )$, 
where $PE := \dot \cup_a P_a E_a$, $\pi'$ is the map from 
$PE$ onto $\Po$, and 
$J'_b v :=$ $J_b v$, 
$v \in P_{\partial_1 b} E_{\partial_1 b}$, $b \in \Si_1(\Po)$, 
$a \in \Si_0(\Po)$.
Since $J_b P_{\partial_1 b} = P_{\partial_0b} J_b$, we conclude that
$
J'_b P_{\partial_1b} E_{\partial_1b} 
\subseteq 
P_{\partial_0b} E_{\partial_0b}$,
so that $J'$ is a well-defined net structure, and $P \Ecal$ is a
Hilbert quasinet bundle.
\end{proof}

Hilbert quasinet bundles associated with projections as in the previous 
lemma are called {\em direct summands}. A Hilbert net bundle $\Ecal$ is 
said to be {\em irreducible} if it does not admit direct summands different
from ${\bf 0}$ and $\Ecal$.
We summarize the results of the present section in the following 
proposition.

\begin{prop}
The category $qH_{net}(\Po)$ with arrows bounded morphisms is a symmetric
tensor {\it C*}-category with subobjects, direct sums, and identity object 
given by the trivial Hilbert net bundle 
$\iota := ( \Po \times \Co , \pi , j , \Po )$. 
(When $\Po$ is pathwise connected) $\iota$ is simple, i.e. $(\iota,\iota) = \Co$.
\end{prop}

In the next lemma, we establish an equivalence between the presence of nowhere
zero global sections and trivial direct summands.
\begin{lemma}
\label{lem_tds}
Let $\Ecal :=$ $( E, \pi , J , \Po )$ be a Hilbert quasinet bundle.
There is a nowhere zero section $\sigma \in S(\Po;\Ecal)$ if and only if
$\Ecal$ has a trivial direct summand of rank one. 
Thus, an irreducible Hilbert net bundle is
nontrivial if and only if it lacks nowhere zero sections.
\end{lemma}

\begin{proof}
If $\Ecal$ has a trivial direct summand, then it is clear that there
is a nowhere zero section $\sigma : \Po \to E$. 
Conversely, since the net structure $J$ involves isometric maps, if there
is such a section $\sigma$ then up to normalization
we may assume that $\left\| \sigma (\dc) \right\| = 1$,
$\dc \in \Po$. Given the trivial net bundle
$\iota := ( \Po \times \Co , p , j , \Po )$,
we define the map
$I_\sigma ( \dc , z ) :=$ $z \sigma (\dc)$
$\dc \in \Po$, $z \in \Co$.
It is clear that $I_\sigma$ is injective. Since $\sigma (\dc) \in E_\dc$,
$\dc \in \Po$, $I_\sigma$ preserves the fibres. We now verify that
$I_\sigma$ preserves the net structure, i.e. 
$I_\sigma \circ j = J \circ I_\sigma$:
for every $b \in \Si_1(\Po)$, we compute
\[
I_\sigma \circ j_b ( \partial_1b , z ) =
I_\sigma ( \partial_0b , z ) =
z \sigma ( \partial_0b ) =
z J_b \circ \sigma (\partial_1b) =
J_b \circ I_\sigma ( \partial_1b , z )
\ ,
\]
where we used the fact that
$J_b \circ \sigma (\partial_1b) = \sigma (\partial_0b)$.
This proves that $I_\sigma$ is a net bundle morphism.
\end{proof}

%

%\begin{oss}
%The previous result may be interesting from the point of view of algebraic quantum
%field theory, in the following sense. Let us consider a space-time $M$ with nontrivial
%$\pi_1(M)$, and a net $\Fcal$ of
%field algebras defined on a suitable poset of open, contractible subsets of $M$.
%%
%If one considers sectors associated with nontrivial representations
%of $\pi_1(M)$ as in [Brunetti-Ruzzi][], then
%one expects the associated local fields to be local sections of a Hilbert net
%bundle $\Ecal$ embedded in $\Fcal$. Lemma \ref{lem_tds} suggests that if $\Ecal$
%is nontrivial and irreducible, then such local fields do not admit
%an extension to the global field algebra.
%\end{oss}
%

\begin{cor}
Let $d \in \No$ and $\Ecal := ( E , \pi , J , \Po )$ be a rank $d$ Hilbert 
net bundle. If there are $d$ linearly independent sections of $\Ecal$, 
then $\Ecal$ is trivial.
\end{cor}

The above corollary shows that the existence of sufficient {\em global} 
sections means that the Hilbert net bundle is trivial, in contrast
to the usual topological setting. The reason is that, a section of a 
Hilbert net bundle not zero over some $\dc \in \Po$ is nowhere 
zero. Of course, this is not true for ordinary vector bundles. 
The Chern functions introduced in Sec.\ref{s_Cclass} measure the 
obstruction to the existence of global sections.

\subsection{Basic Properties of K-theory.}

In the sequel, we always assume that our Hilbert net bundles are 
finite-dimensional.

Let $\Ecal$ be a Hilbert net bundle over $\Po$. We denote 
the isomorphism class of $\Ecal$ by $\left\{ \Ecal \right\}$. 
Direct sum and tensor product induce a natural
semiring structure on the set $\dot{V}_{net} (\Po)$ of isomorphism classes of 
finite-dimensional Hilbert net bundles. We define $\Krm^0(\Po)$ to be
the Grothendieck ring associated with $\dot{V}_{net} (\Po)$, and denote 
the semiring morphism assigning to each isomorphism class 
$\left\{ \Ecal \right\}$ the associated element of $\Krm^0(\Po)$ by
\begin{equation}
\label{eq_proj_K}
\dot{V}_{net} (\Po) \to \Krm^0 (\Po)
\ \ , \ \
\left\{ \Ecal  \right\} \mapsto [\Ecal]
\ .
\end{equation}
By definition, every element of 
$\Krm^0(\Po)$ can be written as $[\Ecal] - [\Fcal]$, with $\left\{ \Ecal \right\}$,
$\left\{ \Fcal \right\} \in V_{net} (\Po)$. 
In analogy with the usual K-theory, we characterize Hilbert net bundles 
$\Ecal , \Ecal' \to \Po$ with $[\Ecal] = [\Ecal']$: as this
result depends only on the definition of Grothendieck ring, we omit the proof.

\begin{lemma}
Let $\Ecal$, $\Ecal'$ be Hilbert net bundles over $\Po$. Then
$[\Ecal] = [\Ecal'] \in \Krm^0(\Po)$ if and
only if there exists a Hilbert net bundle $\Fcal$ such that $\Ecal \oplus \Fcal$
is isomorphic to $\Ecal' \oplus \Fcal$.
\end{lemma}

%
%\begin{proof}
%If there exists $\Fcal$ as above, then $[\Ecal] + [\Fcal] = [\Ecal'] + [\Fcal]$
%and $[\Ecal] = [\Ecal']$. Vice versa,
%recall that the additive structure of the Grothendieck ring is defined in the
%following way: the
%elements of $\Krm^0(\Po)$ are classes of pairs
%$( \left\{ \Ecal \right\} , \left\{ \Fcal \right\} )$,
%modulo the equivalence relation
%$( \left\{ \Ecal \right\} , \left\{ \Ecal \right\} ) = ( 0 , 0 )$,
%where $0$ is the isomorphism class of the zero net bundle. Thus, if 
%$[\Ecal] = [\Ecal']$
%then
%%
%\[
%( \left\{ \Ecal \right\} , \left\{ \Ecal' \right\} ) =
%( \left\{ \Fcal \right\} , \left\{ \Fcal \right\} )
%\]
%%
%\noindent for some $\left\{ \Fcal \right\} \in V_{net}(\Po)$. We conclude that
%%
%$
%( \left\{ \Ecal \oplus \Fcal \right\} , \left\{ \Ecal' \oplus \Fcal \right\} ) = (0,0)
%$,
%%
%i.e. $\left\{ \Ecal \oplus \Fcal \right\} = \left\{ \Ecal' \oplus \Fcal \right\}$.
%\end{proof}
%

The next two results indicate a drastic difference between ordinary K-theory and 
net K-theory. In fact, we are going to prove that every Hilbert net bundle with 
a complement is trivial (in ordinary K-theory, every vector bundle has a complement,
see \cite[I.6.5]{Kar}).

\begin{prop}
Let $\Tcal_d$ be the trivial rank $d$ Hilbert net bundle. Then every direct
summand of $\Tcal_d$ is trivial.
\end{prop}

\begin{proof}
%
%A first argument for the proof is the following: 
%let $\pi_d : \pi_1(\Po,a) \to \Uo (d)$
%be the representation associated with $\Tcal_d$ in the sense of Prop.\ref{prop_vs_sc}.
%Then $\pi_d$ is trivial. Every direct summand of $\Tcal_d$ is associated with a
%subrepresentation of $\pi_d$. But every subrepresentation of $\pi_d$ is itself trivial.
%We now give a second proof of our proposition, expressed in the language of Hilbert net
%bundles. 
%
Let $\Ecal'$ be a Hilbert net subbundle of $\Tcal_d$, and
$P \in (\Tcal_d,\Tcal_d)$ the projection associated with $\Ecal'$.
Then $P$ is a section of the bundle $B(\Tcal_d,\Tcal_d)$ (see Lemma \ref{lem_sbm}).
Now, $B(\Tcal_d,\Tcal_d)$ is isomorphic to the trivial net bundle
$\Po \times {\mathbb{M}}(d)$ (where ${\mathbb{M}}(d)$ denotes the matrix algebra
of order $d$). Thus, every section of $B(\Tcal_d,\Tcal_d)$ is constant. In particular,
$P$ is a constant section whose range is a fixed vector
subspace $V$ of $\Co^d$, and $\Ecal'$ is isomorphic to the trivial bundle
$\Po \times V$.
\end{proof}

\begin{cor}
Let $\Ecal , \hat \Ecal$ be Hilbert net bundles such that
$\Ecal \oplus \hat \Ecal \simeq \Tcal_n$ for some $n \in \No$. 
Then $\Ecal , \hat \Ecal$ are trivial.
\end{cor}

Hilbert net bundles $\Ecal$, $\Fcal$ (not necessarily of the same rank) are
said to be {\em stably equivalent} if $\Ecal \oplus \Tcal_j$ is isomorphic to
$\Fcal \oplus \Tcal_k$ for some $j,k \in \No$. It follows from the previous
corollary that if a Hilbert net bundle is stably equivalent to a trivial net
bundle, then it is trivial. The set $V^s_{net}(\Po)$ of stable equivalence
classes is an Abelian semiring w.r.t. the operation of direct sum and tensor
product. By definition, there is an epimorphism
\begin{equation}
\label{eq_dv_vs}
\dot{V}_{net}(\Po)  \to  V^s_{net}(\Po)
\end{equation}
mapping each $\{ \Ecal \} \in \dot{V}_{net}(\Po)$ into its stable equivalence
class. The kernel of (\ref{eq_dv_vs}) is given by the set of isomorphism classes
of trivial Hilbert net bundles, hence labelled by the non-negative integers.
For the notion of representation ring, we recommend \cite{Seg} to the
reader; the following result is a direct consequence of \cite[Thm.2.8]{Ruz05}.

\begin{teo}
\label{thm_knet}
Let $a \in \Si_0(\Po)$. 
Then the category $V_{net}(\Po)$ is equivalent to the category of unitary, 
finite-dimensional representations of $\pi_1(\Po,a)$, so that $\Krm^0 (\Po)$ 
is isomorphic to the representation ring of $\pi_1 (\Po,a)$.
\end{teo}

\begin{cor}
\label{cor_dec_irr}
Let $\Ecal :=$ $( E , \pi , J , \Po )$ be a Hilbert net bundle. Then there is a 
unique decomposition $\Ecal =$ $\oplus_{i=1}^r n_i \Vcal_i$, where each $\Vcal_i$
is irreducible, $n_i \in \No$, and $n_i \Vcal_i$ denotes
the direct sum of $n_i$ copies of $\Vcal_i$.
\end{cor}

\begin{proof}
In fact, the above decomposition corresponds to the decomposition of the 
representation $\chi : \pi_1(\Po,a) \to \Uo(d)$ associated with $\Ecal$ 
into irreducibles. Note that a Hilbert net bundle is irreducible if and only if 
the associated representation of $\pi_1(\Po,a)$ is irreducible.
\end{proof}

\begin{cor}
\label{cor_p1a}
Let $( \Po , \leq )$ be a poset such that $\pi_1(\Po,a)$ is Abelian. 
Then every Hilbert net bundle over $\Po$ is a direct sum of line net bundles. 
The ring $\Krm^0(\Po)$ is generated by $\Hrm^1 ( \Po , \To )$ 
as a $\Zo$-module.
\end{cor}

Let us now consider the {\em rank function} $\rho : \dot{V}_{net}(\Po) \to \No$, 
assigning to the (isomorphism class of a) Hilbert net bundle the corresponding rank. 
It is clear that $\rho$ is a semiring epimorphism, so that it makes sense to 
consider the associated extension
$\rho : \Krm^0(\Po) \to \Zo$.
The kernel of $\rho$ is called the {\em reduced}
net K-theory of $\Po$ and denoted by $\widetilde \Krm^0(\Po)$. 
In this way we have a direct sum decomposition
\begin{equation}
\label{eq_ktk}
\Krm^0(\Po) = \Zo \oplus \widetilde \Krm^0(\Po)
\ .
\end{equation}
so that $\widetilde \Krm^0(\Po)$ embeds into $\Krm^0(\Po)$. The reduced group
$\widetilde \Krm^0(\Po)$ encodes the nontrivial part of the K-theory of 
$\Po$: if $\Po$ is simply connected, then $\widetilde \Krm^0(\Po) = 0$.
By the definition of stable equivalence class and recalling 
(\ref{eq_dv_vs},\ref{eq_ktk}),
we conclude that the canonical map $\dot{V}_{net}(\Po) \to \Krm^0(\Po)$ 
factorizes through $V^s_{net}(\Po)$, in such a way that the following 
diagram commutes:
\begin{equation}
\label{eq_dia_k}
\xymatrix{
\dot{V}_{net}(\Po)
\ar[r]
\ar[d]
& V^s_{net}(\Po)
\ar[d]
\\ \Krm^0(\Po)
& \widetilde \Krm^0(\Po)
\ar[l]
}
\end{equation}
Net K-theory satisfies natural functorial properties. 
If $\eta : \Po' \to \Po$ is a poset morphism, 
then the pullback induces a ring morphism 
$\eta^* : \Krm^0(\Po) \to \Krm^0(\Po')$. 
If $\Xcal :=$ $(X,\eta,J,\Po)$ is a poset net bundle with fibre $F$, 
then by (\ref{eq_ex_p}) we find an exact sequence of rings
\begin{equation}
\label{eq_ex_s_k}
\Zo \to
\Krm^0(\Po) \stackrel{\eta^*}{\to}
\Krm^0(X) \to
\Krm^0(F) \ .
\end{equation}
Note that the images of $\Zo$ w.r.t. the above maps correspond to trivial 
representations of the homotopy groups $\pi_1(\Po,a)$, $\pi_1(X,x)$, $\pi_1(F,a')$; 
so that, (\ref{eq_ex_s_k}) restricts to an exact sequence of reduced 
$\Krm$-groups
\[
{\bf 0} \to
\widetilde \Krm^0(\Po) \stackrel{\eta^*}{\to}
\widetilde \Krm^0(X) \to
\widetilde \Krm^0(F) \ .
\]
We conclude that if $\Krm^0(F) = \Zo$ (i.e. $\widetilde \Krm^0(F) = {\bf 0}$), 
then 
$\widetilde \Krm^0(X)$ and $\widetilde \Krm^0(\Po)$ are isomorphic.

\

Let us now consider a rank $d$ Hilbert net bundle $\Ecal :=$ $( E , \pi , J , \Po )$.
For every fibre $E_a$, $a \in \Si_0(\Po)$, we consider the associated projective
space $PE_a$. Each unitary $J_b$, $b \in \Si_1(\Po)$, defines a
homeomorphism $[J]_b : PE_{\partial_1b} \to PE_{\partial_0b}$.
Defining $PE :=$ $\dot{\cup}_a PE_a$ and the obvious
projection $[\pi] : PE \to \Po$, we obtain a net bundle of topological spaces
\begin{equation}
\label{def_pe}
\mathcal{PE} := ( PE , [\pi] , [J] , \Po )
\ ,
\end{equation}
called the {\em projective net bundle} associated with $\Ecal$. We denote by
\[
\mcPE := ( \mPE , [\pi]_\prec , [J] , \Po  )
\]
the poset net bundle associated with $\mathcal{PE}$ in the sense of Sec.\ref{pullback}.
\begin{teo}[The Thom isomorphism]
\label{thm_thom}
Let $( \Ecal , \pi , J , \Po)$ be a rank $d$ Hilbert net bundle. 
For every group $G$, the pullback over $\mPE$ induces 
isomorphisms
\begin{equation}
\label{eq_thom_k}
[\pi]_\prec^* : \Hrm^1 ( \Po,G )
\to
\Hrm^1 ( \mPE,G )
\ \ , \ \
[\pi]_\prec^* : \Krm^0 (\Po)
\to
\Krm^0(\mPE)
\ .
\end{equation}
\end{teo}

\begin{proof}
Since the projective space is pathwise and simply connected, we apply Cor.\ref{cor_CV},
Cor.\ref{cor_top}, and conclude that the map $[\pi]_\prec : \mPE \to \Po$
induces an isomorphism
\begin{equation}
\label{eq_pi1_p}
[\pi]_{\prec,*} : \pi_1 ((\mPE,v) \to \pi_1 (\Po,a)
\ \ , \ \ v \in \mPEa
\ .
\end{equation}
Since we find 
\[
\Zrm^1 ( \Po,G ) \simeq \Hom ( \pi_1(\Po,a),G )
\ \ , \ \
\Zrm^1 ( \mPE,G ) \simeq \Hom ( \pi_1((\mPE,v),G )
\ ,
\]
for every group $G$, the theorem is proved at the level of cohomology. 
The isomorphism at the level of K-theory follows by 
Thm.\ref{thm_knet} and (\ref{eq_pi1_p}).
\end{proof}

In contrast to ordinary geometry, the Thom isomorphism really is an
isomorphism and not a monomorphism. But this result is not as useful as
its topological counterpart, since there is not a well-defined notion of
canonical line {\em net} bundle with base $\mPE$. We shall return to this
point in the sequel (Rem.\ref{ex_ptb}).

\

For every topological space $Y$, we denote the category of
vector bundles over $Y$ by $V_{top}(Y)$, and the subcategory (not full) 
of locally constant vector bundles with arrows locally constant morphisms
by $V_{lc}(Y)$ (see Sec.\ref{sec_lcb} or \cite[Ch.I.2]{Kob}).
Now, a Hilbert net bundle $\Ecal := ( E  , \pi , J ,\Po )$ is 
a net bundle of topological spaces in a natural way, so there is 
an associated vector bundle $\pi : \tau E \to \tau \Po$ 
(see Thm.\ref{teo_net_top}).

\begin{teo}
\label{teo_net_lc}
For every Hilbert net bundle $\Ecal :=$ $( E  , \pi , J ,\Po )$, the projection
$\pi$ defines a continuous map $\pi : {\tau E} \to {\tau \Po}$, and
${\tau E}$ becomes a locally constant vector bundle over ${\tau \Po}$.
Thus, there is a functor
\begin{equation}
\label{def_net_lc}
\tau_* : V_{net}(\Po) \to V_{top}({\tau \Po})
\ \ , \ \
\Ecal \mapsto {\tau E}
\ ,
\end{equation}
providing an isomorphism $V_{net}(\Po) \simeq V_{lc}({\tau \Po})$.
For every net bundle of topological spaces $\Xcal :=$ $( X , p , \Phi , \Po )$, 
the topological pullback $p^* (\tau E)$ defines a functor
\begin{equation}
\label{def_pb_nlc}
p^* : V_{net}(\Po) \to V_{lc}({\tau X})
\ .
\end{equation}
\end{teo}

\begin{proof}
$\pi : {\tau E} \to {\tau \Po}$ is a vector bundle as a 
direct consequence of Thm.\ref{teo_net_top}. The functoriality of the 
map (\ref{def_net_lc}) follows by Lemma \ref{lem_net_tau}.
Applying Cor.\ref{cor_h1_net_top} to the $\Uo(d)$--cocycle associated 
with $\Ecal$ we see that ${\tau E}$ is locally constant. 
Finally, $p^*$ takes values in $V_{lc}(\tau X)$ since 
the pullback of a locally constant bundle is locally constant.
\end{proof}

We emphasize the fact that ${\tau \Po}$ is just a $T_0$-space, thus we
cannot use the usual machinery of differential geometry, as in 
\cite[Ch.I]{Kob} for example, to study properties of the locally 
constant bundle ${\tau E}$.

\begin{oss}
\label{ex_ptb}
The fibre bundle $[\pi] : \tau PE \to \tau \Po$ has simply connected fibres,
thus the topological version of Thm.\ref{teo_exp1}, and the fact that
$\tau PE$ is locally constant, imply that we have isomorphisms
\begin{equation}
\label{eq_TTI}
[\pi]_*  :  \pi_1(\tau PE) \to \pi_1(\tau \Po)
\ \ , \ \
[\pi]^* : V_{lc}(\tau \Po) \to V_{lc}(\tau PE)
\ .
\end{equation}
Now, let us consider the topological pullback $[\pi]^*(\tau E) \to \tau PE$. 
By definition, $[\pi]^*(\tau E)$ has as elements pairs 
$(v,\xi) \in E_a \times PE_a$, $a \in \tau \Po$, and we can define 
the canonical line bundle $\tau L \to \tau PE$ by
\[
\tau L := \{ (v,\xi) \in [\pi]^*(\tau E)  : v \in \xi \}
\ .
\]
As in \cite[Ch.IV.2]{Kar}, the restriction of $\tau L$ over $PE_a$,
$a \in \tau \Po$, is just the usual canonical line bundle
$L_a \to PE_a$. Since $L_a$ is not locally constant
$\tau L$ is not locally constant, so it does not belong 
to the image of $[\pi]^*$, and $[\pi]^*(\tau E) = \tau L \oplus E'$
is direct sum of not locally constant bundles. We conclude that
the splitting principle cannot be applied in the category of
Hilbert net bundles.
\end{oss}

\subsection{Chern classes for Hilbert net bundles.}
\label{s_Cclass}

As we saw in the previous section, the algebraic properties of the category of 
Hilbert net bundles over a poset $\Po$ differ drastically from those of 
vector bundles over a topological space.
This fact is also reflected in the construction of Chern classes.
As we saw in Theorem \ref{teo_net_lc} and sha ll see in \S \ref{sec_lcb}, 
Hilbert net bundles are stricty related to locally constant vector bundles, 
which have trivial Chern classes when the base space is a 
manifold.

A different approach to Chern classes may be based on the picture of Hilbert
net bundles as unitary representations of the homotopy group. From this point
of view, there are some results (\cite{Tho86,Sym91}) on
Chern classes associated with representations of a discrete group $G$
(in our case, the homotopy group).
But unfortunately, such Chern classes do not suit our purpose, 
as they are expressed in terms of the group cohomology of $G$,
or, equivalently, in terms of the singular cohomology of a suitable 
Eilenberg-McLane space associated with $G$. 
In both the cases, an immediate interpretation in terms of properties
of the initial poset is lost; moreover, these classes vanish when the
above-mentioned space is homotopic to a manifold. 

In the following we define analogues of Chern classes,
in terms of homotopy-invariant complex functions on the path
groupoid of the poset.

\subsubsection{The first Chern class.}

Let $d \in \No$ and $\Ecal$ be a Hilbert net bundle
of rank $d$. According to (\ref{eq_b_c}), $\Ecal$ is characterized by
a cohomology class $z \in \Hrm^1(\Po,\Uo(d))$.
By the functoriality of $\Hrm^1(\Po , \ \cdot \ )$, the
determinant map $\det : \Uo (d) \to \To$ induces a map
$\det_* : \Hrm^1 ( \Po , \Uo(d) ) \to \Hrm^1 ( \Po , \To )$.
We define {\em the first Chern class of} $\Ecal$ to be
\[
c_1 (\Ecal) := {\det}_*  z \in \Hrm^1 ( \Po , \To  )
\ .
\]
By the elementary properties of determinants, we find
$c_1 ( \Ecal \oplus \Ecal'  ) = c_1(\Ecal) c_1(\Ecal')$,
and obtain an epimorphism of Abelian groups
\begin{equation}
\label{eq_k_c1}
c_1 : \Krm^0(\Po) \to \Hrm^1 ( \Po , \To )
\ .
\end{equation}
The first Chern class is natural in the sense that, if 
$\eta : \Po' \to \Po$ is a morphism, then
\begin{equation}
\label{eq_f_c1}
\eta^* c_1 (\Ecal) = c_1 (\eta^* \Ecal) \ .
\end{equation}
From a categorical point of view, the first Chern class encodes the cohomological 
obstruction for $\Ecal$ to be a special object: if $c_1(\Ecal) = 0$, then the totally 
antisymmetric line net bundle $\lambda^d \Ecal$ is trivial, and every normalized section
$R \in S(\Po ;\lambda^d \Ecal)$ is a solution of the equation (\cite[(3.19)]{DR89}).

Now, $c_1(\Ecal)$ can be regarded as a morphism from $\pi_1(\Po,a)$ into $\To$.
By Thm.\ref{teo_p1h1}, we conclude that $c_1(\Ecal)$ factorizes through a morphism
\begin{equation}
\label{eq_acc}
\hat c_1 (\Ecal) \in {\mathrm{Hom}} (\Hrm_1(\Po,\Zo),\To)
\ ,
\end{equation}
that we call the {\em Abelianized first Chern class}. By Cor.\ref{cor_h1_p1} 
$\Hrm^1(\Po,\To)$ is isomorphic to ${\mathrm{Hom}} (\Hrm_1(\Po,\Zo),\To)$,
thus it is essentially equivalent to consider $\hat c_1 (\Ecal)$
instead of $c_1(\Ecal)$; but since $\Hrm_1(\Po,\Zo)$ is generally easier to
compute, it is usually convenient to use (\ref{eq_acc}).

\subsubsection{Chern K-classes}
\label{sec_cKc}

We now adapt a classical construction to net K-theory (\cite[IV.2.17]{Kar}).
Let $\Ecal :=$ $(E,\pi,J,\Po)$ be a rank $d$ Hilbert net bundle.
We consider the antisymmetric tensor powers $\lambda^k \Ecal$, 
$k = 1 , \ldots , d$, and define
\begin{equation}
\label{eq_kcc}
{\mathrm{k}}_i(\Ecal)
\ := \
\sum_{k=0}^i
(-1)^k
\binom{d-k}{i-k}
[\lambda^k \Ecal]
\ \ , \ \
i = 1 , \ldots , d
\ .
\end{equation}
Some elementary computations involving the dimensions of antisymmetric 
tensor powers imply that $\rho({\mathrm{k}}_i(\Ecal)) = 0$, thus 
${\mathrm{k}}_i(\Ecal) \in \widetilde \Krm^0(\Po)$.
We call the classes ${\mathrm{k}}_i(\Ecal)$ the {\em Chern K-classes}.
Keeping (\ref{eq_k_c1}) in mind, we find
$c_1 ( {\mathrm{k}}_1 (\Ecal) ) = - c_1 [\Ecal]$.
Applying the well-known identity
\begin{equation}
\label{eq_lilj}
\lambda^i(\Ecal \oplus \Ecal') 
\ = \
\oplus_{l+m=i} \lambda^l \Ecal \otimes \lambda^m \Ecal'
\end{equation}
to $\Ecal' = \Tcal_k$, we see
that if $\Ecal$ has rank $d$ and admits a trivial, rank $k$ direct 
summand, then ${\mathrm{k}}_i (\Ecal) = 0$, $k \leq i \leq d$.
%
%
%\begin{proof}
%It suffices to prove our result for ${\mathrm{k}}_d (\Ecal)$,
%in fact the other cases follow by iteration.
%Let $\Tcal_1$ denote the trivial, rank $1$ Hilbert net bundle, 
%so that our task is to prove that 
%${\mathrm{k}}_{d+1} (\Ecal \oplus \Tcal_1) = 0$.
%The identity (\ref{eq_lilj}) yields
%$\lambda^{d+1}(\Ecal \oplus \Tcal_1) =$ $\lambda^d \Ecal$, and
%$\lambda^k (\Ecal \oplus \Tcal_1) =$
%$\lambda^k \Ecal \oplus \lambda^{k-1} \Ecal$,
%$1 \leq k < d$.
%%
%So that,
%%
%${\mathrm{k}}_{d+1} (\Ecal \oplus \Tcal_1 ) =$
%$1 +
%\sum_{k=1}^d (-1)^k [ \lambda^{k-1}\Ecal \oplus \lambda^k\Ecal ] +
%(-1)^{d+1} [ \lambda^d \Ecal ] = 0$.
%\end{proof}
%
%
We define the {\em total Chern K-class} to be
\[
{\mathrm{k}} (\Ecal) 
:= 
1 + \sum \limits_{i=1}^d {\mathrm{k}}_i(\Ecal) h^i
\ \ , \ \
\Ecal \in V_{net}(\Po)
\ .
\]
By definition, ${\mathrm{k}} (\Ecal) =$
${\mathrm{k}} (\Ecal \oplus \Tcal_k)$ for every $k \in \No$. This fits in 
well with the idea that ${\mathrm{k}} (\Ecal)$ should encode the nontrivial
properties of $\Ecal$, and implies that the classes ${\mathrm{k}}_i$
are well-defined for elements of $V^s_{net}(\Po)$. The above considerations
and the naturality of the pullback yield

\begin{prop}
\label{prop_KC}
We have
\begin{equation}
\label{eq_prop_KC}
{\mathrm{k}}_i ( \Ecal \oplus \Ecal')
\ = \
\sum_{l+m=i}^{d+d'} {\mathrm{k}}_l (\Ecal) \ {\mathrm{k}}_m(\Ecal')
\ \ , \ \
\Ecal , \Ecal' \in V_{net}(\Po)
\ .
\end{equation}
Thus, the total Chern K-class factorizes through a morphism
\[
{\mathrm{k}} : V^s_{net}(\Po) \to 1 + h \widetilde \Krm^0(\Po) [[h]]
\ \ , \ \
{\mathrm{k}} ( \Ecal \oplus \Ecal') =
{\mathrm{k}} ( \Ecal ) \ {\mathrm{k}} ( \Ecal' )
\ ,
\]
such that $\eta^* {\mathrm{k}} (\Ecal) =$
${\mathrm{k}} ( \eta^* \Ecal )$ for every poset morphism
$\eta : \Po' \to \Po$.
\end{prop}

%
%\begin{proof}
%Let $d,d'$ denote respectively the rank of $\Ecal$, $\Ecal'$.
%Applying (\ref{eq_lilj}) we obtain (\ref{eq_prop_KC}), and 
%this implies that the total class ${\mathrm{k}}$ extends to the
%desired morphism.
%%
%Finally, the naturality of ${\mathrm{k}}$ is obvious by the definition
%of pullback of a Hilbert net bundle.
%\end{proof}
%

\subsubsection{Chern functions.}
\label{sec_cf}

Let $\Pi_1 (\Po)$ denote  the set of paths of $\Po$. 
Since each $1$--simplex is a path of length $1$, 
$\widetilde \Si_1(\Po)$ is contained in $\Pi_1(\Po)$.
Given a set $S$, a map $f : \Pi_1 (\Po) \to S$ is said to be
{\em homotopy-invariant} if $f(p) = f(p')$ whenever $p$ is
homotopic to $p'$.

Let $\Ecal :=$ $( E , \pi , J , \Po )$ be a rank $d$ Hilbert net bundle
with associated $\Uo(d)$-cocycle $z \in$ $\Zrm^1 ( \Po , \Uo(d) )$. 
By \cite[Eq.32]{RR06}, $z$ can be extended to a $\Uo(d)$-valued,
homotopy-invariant map on $\Pi_1 (\Po)$; 
in particular, if $p$ is homotopic to a constant
path, then $z(p)$ is the identity $1 \in \Uo(d)$.
Using the trace map ${\mathrm{Tr}}$ and the exterior powers 
$\wedge^k$, $k = 1 , \ldots , d$, we define the maps
\[
\chi^k_z (p) := {\mathrm{Tr}} \wedge^k z(p)
\ \ , \ \
p \in \Pi_1 (\Po)
\ .
\]
Clearly, the restriction of $\chi^k_z$ to $\widetilde \Si_1(\Po)$
yields an element of $\Crm^1 ( \Po , \Co )$.
On the other hand, if we restrict $\chi^k_z$ to $\Po(a)$, $a \in \Si_0(\Po)$,
then by homotopy invariance we find that $\chi^k_z$ can be regarded 
as the character of the representation of
$\pi_1(\Po,a)$ associated with $\lambda^k \Ecal$, $k = 1 , \ldots , d$.
In particular, since the trace is the identity map for rank one 
representations, we find
\begin{equation}
\label{eq_chi_c1}
\chi^d_z = c_1 (\Ecal)
\ \ , \ \
d = \rho (\Ecal)
\ .
\end{equation}
Let us denote the ring of bounded, homotopy-invariant maps 
from $\Pi_1 (\Po)$ to $\Co$, by $\Rrm^1(\Po,\Co)$.
By analogy with the previous section, we introduce the 
{\em Chern functions}
\begin{equation}
\label{eq_taui}
\cf_i \Ecal \in \Rrm^1(\Po,\Co)
\ \ , \ \
\cf_i \Ecal
\ := \
\sum_{k=0}^i
(-1)^k
\binom{d-k}{i-k}
\chi^k_z
\ \ , \ \
i = 1 , \ldots , d
\ .
\end{equation}
For $i > d$, we set $\cf_i \Ecal := 0$.
If $\Tcal_d$ is the trivial Hilbert net bundle of rank $d$, 
then (\ref{eq_lilj}) implies
$\cf_i \Tcal_d = 0$, $i = 1 , \ldots , d$.
If $\Lcal \in$ $V_{net}(\Po)$ is a line net bundle with $\To$-cocycle
$\zeta$, then from (\ref{eq_chi_c1}) we find
\begin{equation}
\label{eq_tau_lb}
\chi^1_\zeta = c_1(\Lcal)
\ \Rightarrow \
\cf_1 \Lcal = 1 - c_1(\Lcal)
\ .
\end{equation}
Since the trace is additive on  direct sums and multiplicative on 
tensor products, the
argument of Prop.\ref{prop_KC} shows that
\begin{equation}
\label{eq_tcc_prod}
\cf_i ( \Ecal \oplus \Ecal' )
\ = \
\sum_{l+m=i} \cf_l \Ecal \cdot \cf_m \Ecal' \ ,
\end{equation}
so that
\begin{equation}
\label{eq_cf}
\cf_{d+k} ( \Ecal \oplus \Tcal_k ) = 0
\ , \
\cf_i ( \Ecal \oplus \Tcal_k ) = \cf_i \Ecal
\ \ , \ \
i,k = 1 , \ldots .
\end{equation}
The previous equalities yield the usual interpretation of the
functions $\cf_i \Ecal$ as obstructions to the triviality of $\Ecal$.
Moreover, (\ref{eq_cf}) implies that each $\cf_i \Ecal$
depends only on the class of $\Ecal$ in $V^s_{net}(\Po)$.
For the first Chern class, we have the following result.

\begin{lemma}
For a Hilbert net bundle $\Ecal$ of rank $d$, we have
\begin{equation}
\label{eq_tau_c1}
c_1(\Ecal)
\ = \
1 + \sum_{i=1}^d (-1)^i \cf_i \Ecal
\ .
\end{equation}
\end{lemma}

\begin{proof}
Let us recall the obvious identity
\begin{equation}
\label{eq_fact}
\sum_{j=0}^N (-1)^j \binom{N}{j} = 0
\ \ , \ \
N \in \No
\ .
\end{equation}
We put $\chi_z^0 := 1$, so that, by the definition of
$\cf_i \Ecal$, the r.h.s. of (\ref{eq_tau_c1}) coincides with
\begin{equation}
\label{eq_fact1}
\sum_{i=0}^d (-1)^i \sum_{k=0}^i (-1)^k \binom{d-k}{i-k} \chi^k_z \ .
\end{equation}
Let us put together the coefficients of each $\chi_z^k$ and set
$j := i-k$, so that in particular $(-1)^{i+k} = (-1)^{j+2k} = (-1)^j$.
Then we find that the quantity (\ref{eq_fact1}) coincides with
\[
\sum_{k=0}^d \left( \sum_{j=0}^{d-k} \binom{d-k}{j} (-1)^j \right) \chi_z^k
\ .
\]
Using (\ref{eq_fact}), we conclude that the terms between the brackets vanishes,
except in the case $d=k$, which provides the coefficient $\chi_z^d$.
Thus, the r.h.s. of (\ref{eq_tau_c1}) coincides with $\chi_z^d$, which
is equal to $c_1(\Ecal)$ by (\ref{eq_chi_c1}).
\end{proof}

\begin{teo}
\label{thm_tCC}
Defining the polynomial
\[
\cf \Ecal (h) 
\ := \
1 + \sum_{i=1}^d \cf_i \Ecal \ h^i
\]
for each $\Ecal \in$ $V_{net}(\Po)$ gives a morphism
\begin{equation}
\label{eq_tcc}
\cf : V^s_{net}(\Po) \to 1 + h \Rrm^1(\Po,\Co)[[h]]
\ \ , \ \
\cf ( \Ecal \oplus \Ecal' ) = {\cf \Ecal} \cdot {\cf \Ecal'}
\ ,
\end{equation}
such that $c_1(\Ecal) = \cf \Ecal (-1)$.
%
%If $\Lcal \in$ $V_{net}(\Po)$ is a line net bundle, then
%%
%$\hat \cf_1 \Lcal =$ $1 - c_1 (\Lcal)$.
%
\end{teo}

\begin{proof}
After Prop.\ref{prop_KC}, (\ref{eq_tcc_prod}) and 
(\ref{eq_tau_c1}), the only nontrivial assertion
that we have to verify is the naturality of (\ref{eq_tcc}).
Let $\eta : \Po' \to \Po$ be a morphism and $z \in \Zrm^1(\Po,\Uo(d))$ the cocycle
associated with $\Ecal$.  Then $\eta^*z \in$ $\Zrm^1(\Po',\Uo(d))$
is the cocycle associated with the pullback $\eta^* \Ecal$, and clearly
$\chi^k_{\eta^*z} = \chi^k_z \circ \eta_1$, where
$\eta_1 : \Pi_1(\Po') \to \Pi_1(\Po)$
is the map induced by $\eta$.
This shows that the functions $\cf_i \Ecal$ are natural, i.e.
$
\cf_i (\eta^* \Ecal) = \eta^* \cf_i \Ecal := \cf_i \Ecal \circ \eta_1
$.
\end{proof}

It is instructive to give the details when $\Ecal =$ $\oplus_i^d \Lcal_i$
is the direct sum of line net bundles. Since $\cf$ is a homomorphism, we get
\begin{equation}
\label{eq_taue_l}
\cf \Ecal (h) \ = \
\prod_{i=1}^d \cf \Lcal_i (h) \ = \
\prod_{i=1}^d \left[ 1 + \hat \cf_1 \Lcal_i h \right] \ = \
\prod_{i=1}^d \left[ 1 + ( 1 - c_1(\Lcal_i) ) h \right]
\ ,
\end{equation}
so that we obtain
\begin{equation}
\label{eq_tau_ab}
\cf_i \Ecal
=
\sum_{1 \leq k_1 < \ldots < k_i \leq d}
\left[ 1 - c_1 (\Lcal_{k_1}) \right]
\cdots
\left[ 1 - c_1 (\Lcal_{k_i}) \right]
\  \ , \ \
i = 1 , \ldots , d \ .
\end{equation}

\begin{oss}
\label{rem_ph1}
Let us regard the Chern functions as maps from
$\pi_1(\Po,a)$ to $\Co$.
Since the trace is invariant under the adjoint action
of $\pi_1(\Po,a)$, each $\cf_i \Ecal$ factorizes through a
map on the orbit space $[ \pi_1(\Po,a) ]$ with elements the classes
$
(p') := \{ p*p'*p^{-1} , p \in \pi_1(\Po,a) \} 
$.
Let $\eta : \pi_1(\Po,a) \to [ \pi_1(\Po,a) ]$ denote the natural projection
and $T : \pi_1(\Po,a) \to \Hrm_1(\Po,\Zo)$ the Abelianization map
defined according to Thm.\ref{teo_p1h1}. 
Since $Tp' = T(p*p'*p^{-1}) \in \Hrm_1(\Po,\Zo)$, $p,p' \in \pi_1(\Po,a)$,
there is a map $\hat T : [ \pi_1(\Po,a) ] \to \Hrm_1(\Po,\Zo)$ such that $T$
can be obtained as the composition
\[
\pi_1(\Po,a) 
\ \stackrel{\eta}{\to} \ 
[ \pi_1(\Po,a) ] 
\ \stackrel{\hat T}{\to} \
\Hrm_1(\Po,\Zo) \ .
\]
Of course, when $\pi_1(\Po,a)$ is Abelian the above maps are isomorphisms
and each $\cf_i \Ecal$ factorizes through a map from $\Hrm_1(\Po,\Zo)$
to $\Co$.
\end{oss}

The following simple result shows how the Chern functions
can be used to compute $V^s_{net}(\Po)$ explicitly when $\Po$ has
homotopy group $\Zo$.
\begin{prop}
\label{prop_p1_Z}
Let $\pi_1(\Po,a) = \Zo$, $a \in \Si_0(\Po)$, and
\[
S := \{ (\chi-1)^{-1} \in \Co : \chi \in \To - \{ 1 \} \}
\ .
\]
Then $V^s_{net}(\Po)$ is isomorphic to the multiplicative 
semigroup of polynomials
\[
P(h) :=
1 + \sum_{i=1}^d a_i h^i
\ \in \
1 + h \Co [[h]]
\ \ , \ \
a_i \in \Co
\ ,
\]
whose zeroes belong to $S$.
\end{prop}

\begin{proof}
Since $\pi_1(\Po,a) = \Hrm_1(\Po,\Zo) =$ $\Zo$ is Abelian, 
every Hilbert net bundle $\Ecal$ is the direct sum of line net bundles
$\Lcal_1 , \ldots , \Lcal _d$ with associated cocycles 
$z_i : \pi_1(\Po,a) \to \To$, $i = 1 , \ldots , d$,
thus  (\ref{eq_taue_l}) applies. 
Now, each $z_i : \Zo \to \To$ is uniquely determined by 
$\chi_i := z_i(1) \in \To$; 
so that, the $d$-ple $\{ \chi_i \}$ $\in \To^d$ is a complete
invariant for $\Ecal$. Eliminating any term $\chi_i = 1$,
we obtain a $r$-ple $\{ \chi_i \} \in$ $\To^r$, $r \leq d$, 
which is a complete invariant for the stable equivalence class of $\Ecal$.
It is now clear that the polynomial 
\[
P(h) := \prod_{i=1}^r \left[ 1 + ( 1 - \chi_i ) h \right]
\ \in \
1 + h \Co [[h]]
\]
has the allowed zeroes. Since the set $\{ \chi_i \}$ 
can be reconstructed from the set of zeroes of $P$, the map
$\{ [\Ecal] \mapsto P \}$ is one-to-one and obviously surjective.
\end{proof}

\subsection{A sketch of equivariant K-theory.}
\label{sec_gKt}

Let $\Po$ be a poset and $G$ a group inducing a left action on $\Po$ by poset 
automorphisms. Then for each $n \in \No$ the set of $n$-simplices 
$\widetilde \Si_n(\Po)$ becomes a $G$-set in the natural way, and we denote the 
images of the action of $g \in G$ on $a \in \Si_0(\Po)$, $b \in \Si_1(\Po)$, 
$c \in \Si_2(\Po)$, by $ga \in \Si_0(\Po)$, $gb \in \Si_1(\Po)$, 
$gc \in \Si_2(\Po)$.

A {\em Hilbert net $G$-bundle} is given by a Hilbert net bundle 
$\Ecal :=$ $( E , p , J , \Po )$ such that: (1) $E$ is endowed with a left $G$-action; 
(2) The projection $p$ is a $G$-map, so that each $g \in G$ induces a bijective map 
$g_a : E_a \to E_{ga}$, $a \in \Si_0(\Po)$; (3) every $g_a$ is a unitary operator 
satisfying the condition
\begin{equation}
\label{eq_gj}
g_{\partial_0 b} \ J_b  \ = \ J_{gb} \ g_{\partial_1 b}
\ , \
b \in \Si_1(\Po)
\ .
\end{equation}
Let $\Ecal$, $\tilde \Ecal$ be Hilbert net $G$-bundles over $\Po$.
A morphism $T \in (\Ecal,\tilde \Ecal)$ is said to be
$G$-{\em equivariant} if $T(gv) = g Tv$, $g \in G$,
$v \in E$.
We denote the set of $G$-equivariant morphisms by $(\Ecal,\tilde \Ecal)_G$
and the associated net bundle by
$B(\Ecal,\tilde \Ecal)_G$.
By construction, $B(\Ecal,\tilde \Ecal)_G$ is a Hilbert net bundle endowed with
a trivial $G$-action.
In this way, we obtain a category $V_{net}(\Po;G)$ naturally endowed
with direct sums and tensor products, and, by using the usual construction,
we can define the {\em equivariant K-theory} $\Krm_G^0(\Po)$.

We now give some elementary properties when $\Po$
has the trivial $G$-action. Firstly, note that
(\ref{eq_gj}) implies that each $g \in G$ defines a section 
of $B(\Ecal ,\Ecal)_G$.
Moreover, when $G$ is compact averaging w.r.t. Haar measure provides 
a projection
\begin{equation}
P_G \in (\Ecal,\Ecal)
\ ,
\end{equation}
and we denote the associated Hilbert net bundle by 
$\Ecal_G := P_G \Ecal \in V_{net}(\Po)$. 
Applying $P_G$ to $\hat \Ecal \otimes \overline \Ecal$, 
we find $P_G B(\Ecal,\hat \Ecal) = B(\Ecal,\hat \Ecal)_G$.
Let $v \in \hat \Ecal$, $w \in \Ecal$,
and $T_{vw} \in B(\Ecal,\hat \Ecal)_G$,  $T_{vw}w' := v \cdot (w,w')$; 
if $(\Ecal,\Ecal)_G \simeq \Co$, then we have a canonical map
\begin{equation}
\label{eq_pg_can}
\Ecal \otimes B(\Ecal,\hat \Ecal)_G 
\to 
\hat \Ecal
\ \ , \ \
w' \otimes T_{vw} \mapsto P_G (T_{w'w}) \cdot v
\ .
\end{equation}

\begin{teo}
\label{teo_KRG}
Let $G$ be a compact group acting trivially on $\Po$. Then there 
is a splitting $\Krm^0_G(\Po) \simeq$ $\Krm^0(\Po) \otimes R(G)$.
\end{teo}

\begin{proof}
We use a classical argument 
(see \cite[Prop.2.2]{Seg} for details):
for each irreducible representation $\alpha$ of $G$,
we consider the associated $G$-module $H_\alpha$
and the corresponding trivial Hilbert net $G$-bundle
$\Tcal_\alpha := ( H_\alpha \times \Po , \pi , j , \Po  )$ 
such that $(\Tcal_\alpha,\Tcal_\alpha)_G \simeq \Co$.
Then we define the map
\begin{equation}
\label{eq_ekt}
R(G) \otimes \Krm^0(\Po) \to \Krm_G^0(\Po)
\ \ , \ \
[\alpha] \otimes [\Ecal] \to [\Tcal_\alpha \otimes \Ecal ]
\ .
\end{equation}
The fact that (\ref{eq_ekt}) is an isomorphism follows by
observing that the isotypical decomposition
\[
\Ecal \simeq
\oplus_\alpha  \   \Tcal_\alpha \otimes B(\Tcal_\alpha,\Ecal)_G
\]
(defined as in (\ref{eq_pg_can})) provides the desired inverse.
\end{proof}

\section{Locally constant bundles and comparison with ordinary $K$-theory.}
\label{sec_lcb}

Let $M$ be a Hausdorff space; we fix a poset $M_{\prec}$ of arcwise and 
simply connected open sets forming a base for the topology of $M$.
A vector bundle $p : E \to M$ is said to be \emph{locally constant} whenever
it admits an atlas 
$$\{ \pi_i : p^{-1}(U_i) \to U_i \times \Co^d \}$$ 
defining \emph{locally constant} transition maps. Given the locally constant vector bundle
$p' : E' \to M$ with atlas $\{ \pi'_{i'} : p'^{-1}(U'_{i'}) \to U'_{i'} \times \Co^{d'} \}$, 
a bundle morphism
$T : E \to E'$
is said to be \emph{locally constant} whenever each map
\[
T_{ii'} : (U_i \cap U'_{i'}) \times \Co^d   \ \to \  (U_i \cap U'_{i'}) \times \Co^{d'}
\ \ , \ \
T_{ii'} := \pi'_{i'} \circ T \circ \pi_i^{-1}
\ ,
\]
is locally of the type $T_{ii'}(x,v) = (x,t_{ii'}v)$, where $t_{ii'} : \Co^d \to \Co^{d'}$
is a \emph{fixed} linear map.

In general a vector bundle on $M$ has \emph{continuous} transition maps, 
and a bundle morphism (also between locally constant bundles) defines \emph{continuous} maps 
\[
t_{ii'} : U_i \cap U'_{i'} \to B(\Co^d , \Co^{d'}) \ ,
\]
so the category $V_{lc}(M)$ of locally constant vector bundles on $M$ 
is a \emph{not} a full subcategory of the category $V(M)$ of vector bundles on $M$.
This is an important point, because locally constant vector bundles
$E$, $E'$ can be isomorphic in $V(M)$ without being isomorphic in $V_{lc}(M)$.
$V_{lc}(M)$ is equipped with
a tensor product and direct sums, so the set $\dot V_{lc}(M)$ of isomorphism
classes of locally constant vector bundles becomes a semiring.

We now recall some well-known facts:

(1) The universal cover $\tilde M$ is a right homogeneous
$\pi_1(M)$-space. For any finite-dimensional unitary representation $u$ 
of $\pi_1(M)$ over a Hilbert space $H_u$, we define 
\begin{equation}
\label{eq_EU}
E_u := 
\tilde M \times_{\pi_1(M)} H_u 
\end{equation}
to be the quotient of $\tilde M \times H_u$ by the equivalence relation 
\[
( x,v ) \sim ( x \gamma , u(\gamma)v  ) 
\ \ , \ \ 
\gamma \in \pi_1(M)
\ .
\]
The projection $p : \tilde M \to M$ induces a projection 
$\pi : E_u \to M$
and it can be verified that $E$ is indeed a locally constant vector bundle.
If $T \in ( u , \hat u )$ is an intertwiner, then the map
\[
T^{lc}(x,v) := (x,Tv) 
\ \ , \ \
x \in \tilde M \ , \ v \in H_u \ ,
\]
induces a morphism from $E_u$ to $E_{\hat u}$.
This construction provides an equivalence from the category of finite-dimensional, unitary
representations of $\pi_1(M)$ onto $V_{lc}(M)$ (see \cite[Sec.1.3]{Hat} or \cite[\S I.2]{Kob}).

(2) Let $M$ be a manifold. A vector bundle $E \to M$ is locally 
constant if and only if it admits a flat connection (\cite[Ch.I]{Kob}).
By (\cite[Prop.II.3.1]{Kob}), this implies that the 
Chern classes of locally constant bundles over manifolds vanish.

The problem to finding non-trivial characteristic classes for flat bundles was solved
by Cheeger and Simons (\cite{CS85}), in the following way.
Given the rank $d$ vector bundle $E \to M$ any Chern class arises from a closed form $c_k(E) \in Z^{2k}_{dR}(M)$ and, 
denoting the group of singular $n$-cycles, $n \in {\mathbb{N}}$, by $Z_n(M)$, there is a unique morphism
\[
c^\uparrow_k(E) : Z_{2k-1}(M) \to {\mathbb{R/Z}}
\ \ , \ \
\varsigma \mapsto \left \langle c^\uparrow_k(E) , \varsigma \right \rangle
\ ,
\]
such that 
\begin{equation}
\label{eq.CS}
\left \langle c^\uparrow_k (E) , \partial \varphi \right \rangle 
\ = \ 
\int_\varphi c_k(E) \ {\mathrm{mod}} \, {\mathbb{Z}}
\ \ , \ \
\forall k = 1 , \ldots , d \ , \ \varphi \in C_{2k}(M) \ ,
\end{equation}
where $C_{2k}(M)$ is the set of singular $2k$-chains and
$\partial : C_{2k}(M) \to Z_{2k-1}(M)$
is the boundary.
When $E$ is locally constant $c_k(E) = 0$, so $c^\uparrow_k(E)$ vanishes on $\partial C_{2k}(M)$ by (\ref{eq.CS}); 
this implies that $c^\uparrow_k(E)$ is a cocycle in $Z^{2k-1}(M,{\mathbb{R/Z}})$,
whose cohomology class is called the {\em secondary Cheeger-Chern-Simons $k$-class} of $E$
(\cite[Theorem 1.1]{CS85}).

Given the space $M$, the odd cohomology with coefficients 
in the ring $R$ is given by the direct sum
\[
H^{odd}(M,R) \ := \ {\mathbb{Z}} \oplus \bigoplus_{k = 1}^\infty H^{2k-1}(M,R) \ .
\]
When $M$ is a manifold $H^i(M,R)$ is eventually trivial and the above direct sum is indeed finite.

\begin{teo}
\label{teo_nb_lcb}
Let $M$ be a locally compact Hausdorff space, locally arcwise and simply connected.
Then the categories $V_{net}(M_{\prec})$ and $V_{lc}(M)$ are equivalent.
Moreover, $\Krm^0(M_\prec)$ is isomorphic to the
representation ring of $\pi_1(M)$, and describes the Grothendieck ring 
associated with $\dot V_{lc}(M)$.
There is a forgetful map into topological K-theory
\begin{equation}
\label{eq.k}
\Krm^0 (M_\prec)  \to  \Krm^0 (M)
\end{equation}
and, when $M$ is a manifold, we have the (additive) Cheeger-Chern-Simons character
\begin{equation}
\label{eq.ccs}
ccs : \Krm^0 (M_\prec)  \to  H^{odd}(M,{\mathbb{R/Q}})
\ .
\end{equation}
\end{teo}

\begin{proof}
As we saw the map $\{ u \mapsto E_u \}$ defined by (\ref{eq_EU})
provides an equivalence between the category $Rep(\pi_1(M))$ of 
finite-dimensional, unitary representations of $\pi_1(M)$ and $V_{lc}(M)$. 
Moreover, $Rep(\pi_1(M))$ is equivalent to $V_{net}(M_\prec)$ by \cite[Prop.3.8]{RR06}.
Thus, $V_{net}(M_\prec)$ is equivalent to $V_{lc}(M)$, and Theorem \ref{thm_knet} implies 
that $\Krm^0(M_\prec)$ is the 
(finite dimensional) representation ring of $\pi_1(M)$.
Moreover, the map (\ref{eq.k}) is defined by embedding $V_{lc}(M)$ into $V(M)$.
Finally (\ref{eq.ccs}) is defined by
\[
ccs ([\Ecal]) \ := \ rank(\Ecal) + \sum_{k=1}^d \frac{ (-1)^{k-1}  }{ (k-1)! } \ [ c^\uparrow_k(E) ]_{ {\mathrm{mod}}  {\mathbb{Q}} }
\ \ , \ \
\forall \Ecal \in V_{net}(M_{\prec})
\ ,
\]
as in \cite[Theorem 8.22]{CS85}, where $E$ is the locally constant vector bundle defined by $\Ecal$.
\end{proof}

The functor 
$F : V_{net}(M_{\prec}) \to V_{lc}(M)$ 
giving the above equivalence can be explicitly described as follows.
If $\Ecal \in V_{net}(M_\prec)$, 
then by Prop.\ref{prop_vs_sc} there is $d \in \No$ 
and a unitary representation 
$z : \pi_1(M_\prec,a) \to \Uo (d)$, 
$a \in \Si_0(M_\prec)$,
defining a representation $u$ of $\pi_1(M)$
according to \cite[Thm.2.18]{Ruz05}
(in concrete terms, if $\gamma : [0,1] \to M$ is a closed curve,
define $u (\gamma) := z (p)$, 
where $p \in M_\prec(a)$ is an {\em approximation}
of $\gamma$ in the sense of \cite[Def.2.13]{Ruz05}).
We define $F(\Ecal) := E_u$.
If $T \in ( \Ecal , \hat \Ecal )$ is a morphism, 
then according the proof of \cite[Prop.3.8]{RR06}
we have an intertwiner
$T(a) \in (z,\hat z) = (u,\hat u)$,
providing a morphism from $E_u$ to $E_{\hat u}$.

\begin{cor}
\label{cor_p1_ab}
Let $\pi_1(M)$ be Abelian. Then
every locally constant vector bundle over $M $decomposes into a direct sum of
locally constant line bundles. The ring $\Krm^0 (M_\prec)$ is generated as a
$\Zo$-module by the singular cohomology $\Hrm^1 ( M ,\To  )$.
\end{cor}

\begin{proof}
Apply Thm.\ref{teo_ho_co}.
\end{proof}

\begin{ex}[The circle]
We discuss the case $M = S^1$ endowed with the base $S^1_\prec$ of open, non-dense intervals.
It is well known that 
\[
\pi_1(S^1) \ \simeq \ \Hrm_1 (S^1,\Zo) \ \simeq \ \Zo \ ,
\]
so we obtain
$\pi_1(S^1_\prec) \ \simeq \ \Zo$
and
\[
\Hrm^1 (S^1,\To) \ \simeq \ \Hrm^1 (S^1_\prec,\To) \ \simeq \ \To \ .
\]
Moreover any locally constant vector bundle is direct sum of locally constant line bundles,
and the same is true for Hilbert net bundles over $S^1_\prec$.
By Cor.\ref{cor_p1a} we have 
\[
\Krm^0(S^1_\prec) \ \simeq \ \Zo [\To] \ ,
\]
where $\Zo [\To]$ is the ring of formal linear combinations of elements of $\To$ with coefficients in $\Zo$.
Since any vector bundle on $S^1$ is trivial we have $\Krm^0 (S^1) \simeq \Zo$, and (\ref{eq.k})
takes the form
\[
\Zo [\To] \to \Zo \ \ , \ \ p \mapsto \sum_{z \in S} n_z \ \ , \ \  \forall p := \sum_{z \in S}^{S \subset \To , |S| < \infty } n_z z \ ,
\]
corresponding to the operation of assigning the rank to the given Hilbert net bundle.
Finally, we have 
\[
H^k(S^1,{\mathbb{R/Q}})
\ \simeq \
\left\{
\begin{array}{ll}
{\mathbb{R/Q}} \ \ , \ \ k = 1 \ , \ 
\\
\{ 0 \}  \ \ , \ \ k > 1 \ ,
\end{array}
\right.
\]
so $H^{odd}(S^1,{\mathbb{R/Q}}) \simeq \Zo \oplus {\mathbb{R/Q}}$ and (\ref{eq.ccs}) takes the form
\[
\Zo [\To] \to \Zo \oplus {\mathbb{R/Q}}
\ \ , \ \
ccs(p)  = 
\left( \sum_{z \in S} n_z \right) \oplus 
\left( \sum_{z \in S} n_z \log z_{ {\mathrm{mod}} {\mathbb{Q}} } \right)
\ \ , \ \ 
\forall p \in \Zo [\To]
\ ,
\]
having used the logarithm $\log : \To \to {\mathbb{R/Z}}$ and the quotient
${\mathrm{mod}} {\mathbb{Q}} : {\mathbb{R/Z}} \to {\mathbb{R/Q}}$.
\end{ex}

The {\em homotopy groupoid of} $M$ is given by the category  
with objects the points of $M$ and set of arrows the set 
$\tilde{\pi}_1(M)$ of homotopy classes of continuous curves in $M$.
We denote the ring of bounded, 
complex functions on $\tilde{\pi}_1(M)$ by $\Rrm^1(M,\Co)$
and 
the Abelian semigroup of stable
equivalence classes of locally constant vector bundles, 
defined as in (\ref{eq_dv_vs}), by $V_{lc}^s(M)$.
\begin{teo}
For each locally constant vector bundle $E \to M$ of rank $d$, there are
functions
$\cf_i E \in \Rrm^1(M,\Co)$, $i = 1 , \ldots , d$,
such that
$\cf_i (E \oplus E') = \sum_{l+m=i} \cf_l E \cdot \cf_m E'$.
The polynomial
\[
\cf E (h) := 1 + \sum_i^d \cf_i E \ h^i
\]
defines a morphism
\[
\cf : V_{lc}^s(M) \to 1 + h \Rrm^1(M,\Co)[[h]]
\ \ , \ \
\cf (E \oplus E') = \cf E \cdot \cf E'
\ .
\]
When $\pi_1(M)$ is Abelian, $\Rrm^1(M,\Co)$ can be replaced by 
the ring of bounded complex functions on the singular homology 
$\Hrm_1(M,\Zo)$.
\end{teo}

\begin{proof}
By \cite[Thm.2.18]{Ruz05}, the homotopy groupoid of $M$ is isomorphic 
to the homotopy groupoid of $M_\prec$, so that there is a ring isomorphism
$\Rrm^1(M_\prec,\Co) \simeq \Rrm^1(M,\Co)$.
Thus, we apply Thm.\ref{teo_ho_co}, Thm.\ref{thm_tCC} and 
Rem.\ref{rem_ph1}.
\end{proof}

%*************************************************************

\appendix

\section{Simplicial sets.}
\label{sec_ssets}

A simplicial set is a contravariant 
functor from the simplicial category $\Delta^+$ to the category of sets. 
$\Delta^+$ is a subcategory of the category of sets having as objects 
$n:=\{0,1,\dots,n-1\}$, $n\in\mathbb N$ and as mappings the order 
preserving mappings. A simplicial set has a well known description 
in terms of generators, the face and degeneracy maps, and relations. 
We use the standard notation $\partial_i$ and $\sigma_j$ for the face 
and degeneracy maps, and denote the compositions
$\partial_{i}\partial_{j}$, $\si_{i}\si_{j}$, respectively, by 
$\partial_{ij}$, $\si_{ij}$. A path in a simplicial set 
is an expression of the form 
\[
p:=b_n*b_{n-1}*\cdots*b_1 \ ,
\]
where the $b_i$ are $1$--simplices 
and $\partial_0b_i=\partial_1b_{i+1}$ for $i=1,2,\dots,n-1$. We set 
$\partial_1p:=\partial_1b_1$, $\partial_0p:=\partial_0b_n$ and 
\[
\ell(p):=n
\]
the length of $p$. Concatenation 
gives us an obvious associative composition law for paths and in this way 
we get a category without units.\smallskip 

Homotopy provides us with an equivalence relation $\sim$ on this structure. 
This is the equivalence relation generated by a finite sequence $s(i)$,
with $i=1,\ldots,k$  say, 
of elementary deformations. An \emph{elementary deformation} of a path 
consists in replacing a $1$--simplex $\partial_1c$ of the path
by a pair $\partial_0c,\partial_2c$, where $c\in\Si_2$, or, conversely 
in replacing a consecutive pair $\partial_0c,\partial_2c$ of $1$--simplices 
of $p$ by a single $1$--simplex $\partial_1c$.  
The former type of deformation is called an
\emph{ampliation} of the path,  the latter  a \emph{contraction}. 
Quotienting by  this equivalence relation yields the 
homotopy category of the simplicial set.\smallskip

We shall mainly use symmetric simplicial sets. 
These are contravariant functors 
from $\Delta^s$ to the category of sets, 
where $\Delta^s$ is the full subcategory 
of the category of sets with the same objects as $\Delta^+$. A symmetric 
simplicial set also has a description in terms of generators and relations, 
 where the generators now include the permutations of adjacent vertices, 
denoted $\tau_i$. In a symmetric simplicial set we define the reverse of a 
$1$--simplex $b$ to be the $1$--simplex $\overline{b}:=\tau_0b$ 
and the reverse of a path $p=b_n*b_{n-1}*\cdots*b_1$ is the path 
$\overline{p}:=\tau_0b_1*\tau_0b_2*\cdots*\tau_0b_n$. The reverse acts as 
an inverse after taking equivalence classes so the homotopy category 
becomes a homotopy groupoid. Given a symmetric simplicial set 
$\widetilde \Si_*$ we denote its homotopy 
groupoid by $\pi_1(\widetilde \Si_*)$.
%******************

\subsection{Homotopy of products}
\label{sec_hom_prod}

Consider a pair $\widetilde{\Si}^\alpha_*$ and 
$\widetilde{\Si}^\mu_*$ 
of symmetric simplicial sets. Let $\Pi^\alpha_1$ and 
$\Pi^\mu_1$ denote the corresponding set of paths, and 
let $\sim_\alpha$ and 
$\sim_\mu$ denote the corresponding 
homotopy equivalence relations. Now, consider the product simplicial set 
$\widetilde{\Si}^\alpha_*\times \widetilde{\Si}^\mu_*$. Let 
$\Pi^{\alpha\times\mu}_1$
be the set of paths of the product simplicial set and denote the 
homotopy equivalence relation by $\sim$. A path 
$p$ in this set is a pair $(p^\alpha,p^\mu)$, where  
$p^\alpha,p^\mu$ are paths in $\Pi^\alpha_1$ 
and $\Pi^\mu_1$ respectively with
$\ell(p^\alpha)=\ell(p^\mu)$.  Note that a homotopy 
in this set is a finite sequence $s(i)=(s^\alpha(i),s^\mu(i))$, with 
$i=1,\ldots,m$ to say,  
where $s^\alpha(i)$ and $s^\mu(i)$ are either 
both ampliations or both contractions 
of the paths $p^\alpha$ and $p^\mu$ for any $i=1,\ldots, m$.
In particular note that, if $p\sim q$, then $p^\alpha\sim_\alpha q^\alpha$ and 
$p^\mu\sim_\mu q^\mu$. \smallskip 

Our aim is to show that the fundamental 
groupoid of the product simplicial set 
is equal to the product of the fundamental groupoids. 
To this end,  consider the set
$\Pi^\alpha_1\times\Pi^\mu_1$ product of paths.
Note that the elements $p$ of this set are pairs 
$(p^\alpha,q^\mu)$  of paths where, in general,  
$\ell(p^\alpha)\ne \ell(p^\mu)$. Finally, 
observe that the proof of our claim follows once we have shown
that  the identity map 
\[
\Pi_1^{\alpha\times \mu} \ni p\to p \in \Pi^\alpha_1\times \Pi^\mu_1 \  , 
\]
is injective and surjective  up to equivalence. It is easily seen to be
surjective. In fact let  $p\in \Pi^\alpha_1\times \Pi^\mu_1$. Assume 
that $\ell(p^\alpha)=\ell(p^\mu) + k$. If $a$ is the starting point 
of the path $p^\mu$ in $\Si^\mu_*$, then we consider the path 
$p^\mu*(\sigma_0a)^{*_k}$, where  $(\sigma_0a)^{*_k}$ is the 
$k$-fold composition of the degenerate 1--simplex $\si_0a$. Clearly 
$p^\mu*(\sigma_0a)^{*_k}\sim_\mu p^\mu$ and 
$(p^\alpha,p^\mu*(\sigma_0a)^{*_k})\in\Pi_1^{\alpha\times\mu}$.\smallskip

For injectivity we must show that 
\begin{equation}
\label{eq_pss}
p^\alpha\sim_\alpha q^\alpha \ \ , \ \ p^\mu\sim_\mu q^\mu
\ \ \Rightarrow \ \
p \sim q
\ .
\end{equation}
If $p'\sim p$ and $q'\sim q$ then it obviously suffices to pose the
question for $p'$ and $q'$ instead. 
$p*\sigma_0\partial_1p$ is homotopic to $p$ so it suffices to 
suppose that $\ell(p)=\ell(q)$. Pick 
sequences $s^\alpha_0$ and $s^\mu_0$ of elementary deformations leading from 
$p^\alpha$ to 
$q^\alpha$ and from $p^\mu$ to $q^\mu$. Since $\ell(p^\alpha)=\ell(q^\alpha)$ 
and $\ell(p^\mu)=\ell(q^\mu)$, each sequence contains as many ampliations as 
contractions but they may not have the same length. If $\ell(s^\alpha_0)<\ell(s^\mu_0)$, 
say, we may adjoin pairs consisting of an ampliation adding on $\sigma_0\partial_1p^\alpha$ 
and a contraction removing it again. We may therefore 
replace $s_0^\alpha$ and $s_0^\mu$ by sequences $s_1^\alpha$ and
$s_1^\mu$ of the same length. We now proceed inductively: if after
$n-1$ steps we have sequences $s_n^\alpha$ and $s_n^\mu$ and if the
$n-th$ elements
of $s_n^\alpha$ and $s_n^\mu$ are both ampliations or both 
contractions, set $s_{n+1}^\alpha:=s_n^\alpha$ and $s_{n+1}^\mu:=s_n^\mu$. If this 
is not the case 
and the $n-th$ element of $s_n^\alpha$, say, is a contraction, $s_{n+1}^\alpha$ is 
obtained by inserting $\sigma_0\partial_1p^\alpha$ as $n-th$ member and $s_{n+1}^\mu$ 
by adding $\sigma_0\partial_1p^\alpha$ as a final element. 
The process terminates with sequences $s_k^\alpha$ and $s_k^\mu$, say, 
having ampliations and contractions in the same relative positions. 
Setting $s(i):=(s_k^\alpha(i),s_k^\mu(i))$. 
$s$ is then a sequence of deformations from $p$ to 
$q*\sigma_0a*\cdots*\sigma_0a\sim q$, 
where $a=\sigma_0\partial_1(p^\alpha,p^\mu)$, completing the proof.

\subsection{Simplicial sets of a poset}
\label{sec_ss_poset}

We shall be concerned here with two different simplicial sets that can be 
associated with a poset and we give their definitions not just for a poset 
but for an arbitrary category $\Ccal$. The first denoted $\Si_*(\Ccal)$ is 
just the usual nerve of the category. Thus the $0$--simplices are just the 
 objects of $\Ccal$, the $1$--simplices are the arrows of $\Ccal$ and a 
 $2$--simplex $c$ is made up of its three faces which are arrows satisfying 
 $\partial_0c\,\partial_2c=\partial_1c$. The explicit form of higher simplices 
will not be needed in this paper. The homotopy category of $\Si_*(\Ccal)$ 
is canonically isomorphic to $\Ccal$ itself. The second simplicial 
set $\tilde\Si_*(\Ccal)$ is a symmetric simplicial set. 
It is constructed as follows. 
Consider the poset $P_n$ of non-void subsets 
of $\{0,1,\dots,n-1\}$ ordered under inclusion. Any mapping $f$ from 
$\{0,1,\dots,m-1\}$ to $\{0,1,\dots,n-1\}$ induces an order preserving 
mapping from $P_m$ to $P_n$. Regarding the $P_n$ as categories, 
we have realized $\Delta^s$ as a subcategory of the category of 
categories. We then get a symmetric simplicial set where an 
$n$--simplex of $\tilde\Si_*(\Ccal)$ is a functor from $P_n$ to
$\Ccal$. \smallskip
Given a poset $\Po$ with order relation $\leq$. We recall that 
$\Po$ is \emph{upward} directed whenever 
for any pair $o,\hat o\in K$ there is $\tilde o$ such that 
$o,\hat o\leq \tilde o$. It is \emph{downward} directed if the dual poset 
$K^\circ$ is upward directed. The \emph{dual} $K^\circ$ of $\Po$ 
is the poset having the same elements as $\Po$ and opposite order relation 
$\leq^\circ$, i.e., $a\leq^\circ \tilde a$ if, and only if, 
$a\geq \tilde a$.\smallskip

When we specialize to a poset $\Po$  the nerve $\Si_*(K)$ and 
$\widetilde\Si_*(K)$ admit the following representation
(see \cite{RRV07} for details).  
A $0$--simplex of $\widetilde\Si_*(K)$ 
is just an element of the poset. 
For $n\geq 1$, an $n-$simplex $x$ is formed
by $n+1$  $(n-1)-$simplices $\partial_0x, \ldots,\partial_nx$,
and by a $0$--simplex $|x|$ called the \emph{support}
of $x$ such that $|\partial_0x|, \ldots,|\partial_{n}x|\leq |x|$.
The nerve $\Si_*(K)$ turns out to be a subsimplicial set of $\tilde\Si_*(K)$. 
To see this, it is enough to define 
$f_0(a):= a$ on  $0$--simplices and, inductively,  
$|f_n(x)|:= \partial_{01\cdots (n-1)}x$ and 
$\partial_i f_n(x):= f_{n-1}(\partial_ix)$. 
So we obtain a simplicial map 
$f_*:\Si_*(K)\to \tilde \Si_*(K)$.
We sometimes adopt the following 
notation: $(o;a,\tilde a)$ is the
1--simplex of $\tilde\Si_1(K)$ whose support is $o$ and 
whose 0-- and 1--face are, respectively, $a$ and $\tilde a$;
$(a,\tilde a)$ is the 1--simplex of the nerve $\Si_1(K)$
whose 0-- and 1--face are, respectively, 
$a$ and $\tilde a$.\smallskip

The poset $\Po$ is said to be pathwise connected 
whenever the simplicial set $\tilde\Si_*(K)$ is pathwise 
connected. The homotopy groupoid $\pi_1(K)$  
of $\Po$ is defined as $\pi_1(\tilde\Si_*(K))$; 
in particular, choosing a reference point $a \in \Si_0(\Po)$
we get the homotopy group $\pi_1(\Po,a)$. We observe that 
$\Po$ is pathwise connected if, and only if, 
its dual $K^\circ$ is pathwise connected and that 
$\pi_1(K)$ is isomorphic to $\pi_1(K^\circ)$ (see \cite{RRV07}). 
In the present paper we will consider only pathwise 
connected posets $\Po$. Thus we shall say that $\Po$ is simply connected whenever 
$\pi_1(K)$ is trivial. Note that, when $\Po$ is upward directed, 
$\widetilde \Si_*(K)$ admits a contracting homotopy. So in this case
$\Po$ is simply connected. The same happens when 
$\Po$ downward directed since  $\Po$ and $K^\circ$ have isomorphic homotopy 
groupoids.

%*************************************************************

%{\small

\end{document}